\documentclass[12pt,reqno]{amsart}
\newcommand{\mysection}[1]{\section{#1}
      \setcounter{equation}{0}}

\usepackage{color}

\newtheorem{theorem}{Theorem}[section]
\newtheorem{lemma}[theorem]{Lemma}
\newtheorem{proposition}[theorem]{Proposition}

\theoremstyle{definition}
\newtheorem{assumption}{Assumption}[section]
\newtheorem{definition}{Definition}[section]

\theoremstyle{remark}
\newtheorem{remark}{Remark}[section]

\newcommand\bR{\mathbb{R}}
\newcommand\bB{\mathbb{B}}
\newcommand\bS{\mathbb{S}}

\newcommand\cI{\mathcal{I}}
\newcommand\cJ{\mathcal{J}}

\newcommand\cR{\mathcal{R}}
\newcommand\cQ{\mathcal{Q}}

\newcommand\cM{\mathcal{M}}
\newcommand\cK{\mathcal{K}}
\newcommand\cN{\mathcal{N}}

\makeatletter
 \newcommand{\sumstar}%
 {\operatornamewithlimits{\sum@\kern-.2em\raise1ex\hbox{*}}}
 \makeatother

\renewcommand{\boldsymbol}[1]{\mathbf{#1}}

\begin{document}

\title[Rate of convergence]{
 An accelerated splitting-up method for parabolic equations}

\author[I. Gy\"ongy]{Istv\'an Gy\"ongy}
\address{School of Mathematics,
University of Edinburgh,
King's  Buildings,
Edinburgh, EH9 3JZ, United Kingdom}
\email{gyongy@maths.ed.ac.uk}

\author[N.  Krylov]{Nicolai Krylov}%
\thanks{The work of the second author was partially supported
by NSF Grant DMS-0140405}
\address{127 Vincent Hall, University of Minnesota,
Minneapolis,
       MN, 55455, USA}
\email{krylov@math.umn.edu}


\begin{abstract}
We 
approximate the solution $u$
of the Cauchy problem
$$
\frac{\partial}{\partial t} u(t,x)=Lu(t,x)+f(t,x),
\quad (t,x)\in(0,T]\times\bR^d,
$$
$$
u(0,x)=u_0(x),\quad x\in\bR^d
$$
by splitting the equation into the system 
$$
\frac{\partial}{\partial t} v_r(t,x)=L_rv_r(t,x)+f_r(t,x),
\qquad r=1,2,...,d_1,
$$
where $L,L_r$ are second order
differential operators, $f$, $f_r$
are functions of $t,x$, such that
$L=\sum_r L_r$, $f=\sum_r f_r$.
Under natural conditions on solvability 
in the Sobolev spaces $W^m_p$, 
we show that for any
$k>1$ one can approximate the solution $u$ with an 
error of order $\delta^k$,  
by an appropriate combination of the solutions $v_r$
along a sequence of time discretization, 
where
$\delta$ is proportional to the step  size of the grid. 
This result is obtained by using
the 
time change introduced in \cite{KG1}, together with 
Richardson's method and 
a power series expansion of the error of 
splitting-up approximations in terms of~$\delta$. 
\end{abstract}

\maketitle

\mysection{Introduction}                   \label{section16.03.09.03}

In this paper we are interested 
in the rate of convergence of  
splitting-up approximations 
to the  
solution of the  parabolic, possibly {\em 
degenerate,\/}
differential equation with time dependent coefficients 
\begin{equation}                          \label{17.03.09.03}
\frac{\partial}{\partial t}
u(t,x)=Lu(t,x)+f(t,x),
\quad (t,x)\in(0,T]\times\bR^d
\end{equation}
with initial condition 
\begin{equation}                          \label{18.03.09.03}
u(0,x)=u_0(x),\quad x\in\bR^d,
\end{equation}
where $L$ is a differential 
operator of the form 
$$
L=a^{ij}(t,x)\frac{\partial^2}
{\partial x^ix^j} + 
a^{i}(t,x)\frac{\partial }{\partial x^i}
+a(t,x)
$$
and $f$ is a function of 
$t\geq 0$ and $x\in\bR^d$. 
The first step in the splitting methods 
is to choose suitable 
decompositions 
$L=L_1+L_2+...+L_{d_1}$ and 
$f=f_1+f_2+...+f_{d_1}$ 
for the operator $L$ and the free term 
$f$, such that each  
equation                                    
\begin{equation}                                \label{19.03.09.03}
\frac{\partial}{\partial t} v_r(t,x)
=L_rv_r(t,x)+f_r(t,x)
\end{equation}
$r=1,2,...,d_1$ is integrable exactly 
or approximately. 
   
Assume for simplicity 
that the operators 
$L_r$ and  the  free terms $f_r$,  
$r=1,2,...,d_1$, 
are time independent, and, 
for fixed $T>0$ and integer $n\geq 1$, 
consider the uniform step   
\begin{equation}                                 \label{grid}
T_n:=\{t_i^n:=iT/n , i=0,1,2,...,n\}
\end{equation} 
of step  size 
$\delta:=T/n$. 
Then a splitting-up
approximation $u^{(n)}$ for the solution $u$ of 
(\ref{17.03.09.03})-(\ref{18.03.09.03}) 
is defined by  
\begin{equation}                              \label{12.05.09.03}
u_{n}(t_i^n)
:=(\bS_{\delta}^{(d_1)}...
\bS_{\delta}^{(2)}\bS^{(1)}_{\delta})^iu_0, 
\qquad i=0,1,...,n
\end{equation}              
at the grid points.  
Here $\bS^{(r)}_{t}$ denotes 
the solution operator of equation 
(\ref{19.03.09.03}), i.e., 
$\bS^{(r)}_{t}\varphi$ is the solution of
(\ref{19.03.09.03})  at time $t$  
with initial condition $\varphi$ at $t=0$. 
Formula (\ref{12.05.09.03}) 
means that we take $u^{(n)}(0)=u_0$, 
and we calculate 
the approximation at a grid point
$t+\delta$ from the approximation 
$u^{(n)}(t)$ at the previous grid point $t$,  
by solving equations (\ref{19.03.09.03}) 
for $r=1,2,...,d_1$ on
the  same time interval $[0,\delta]$  
successively. First we solve the first equation ($r=1$) 
on $[0,\delta]$ with initial condition  
$v_1(0)=u^{(n)}(t)$, 
and then we solve the 
second equation, third equation, and so on, 
on the same interval 
$[0,\delta]$, by taking 
always the value at $\delta$ of 
the solution of the previous
equation as the initial value for 
the following equation. 
Finally we solve the last 
equation ($r=d_1$) on the interval
$[0,\delta]$ with initial condition 
$v_{d_1}(0)=v_{d_1-1}(\delta)$, and the 
value $v_{d_1}(\delta)$  
is the value of the
splitting-up approximation at $t+\delta$. 

This kind of approximations is well-known in numerical
analysis, and it has been successfully applied to various
types of PDE problems. 
They are often combined with other numerical
methods, as finite differences, finite elements, etc. 
Pioneering applications to the heat equation, 
to hyperbolic equations,
to nonlinear PDEs are presented, 
for example, in
\cite{PR}, 
\cite{DR}, \cite{BG}, 
in \cite{Ya}, \cite{HT} and in \cite{C}, \cite{BM}
\cite{T},  
respectively. 
Many applications and modifications of 
the splitting-up method have
been developed in various applied fields 
of linear and nonlinear 
PDEs and ODEs, under 
a variety of different names, 
like dimensional splitting, 
operator splitting, 
predictor-corrector method, 
method of alternating  directions, 
fractional step method, 
Lie-Trotter-Kato formula, Baker-Campbell-Hausdorff formula, 
Chernoff formula, split Hamiltonian,  split-steps, leapfrog. 
For guidance in the huge varieties of methods, 
names and references 
we refer to the survey article 
\cite{MQ} and books \cite{M1}, \cite{M2}.  

In the context of semigroups the splitting-up 
method first appears
as Trotter's formula \cite{T}, 
which can be formulated as follows: 
$$    
\lim_{n\to\infty}(e^{tA_{d_1}/n}
...e^{tA_2/n}e^{tA_1/n})^nz=e^{tA}z, 
\quad \forall z\in\bB, 
$$
where $A=A_1+A_2+...+A_{d_1}$ and $A_r$ 
are infinitesimal generators 
of $C_0$-semigroups of contractions 
$\{e^{tA}:t\geq 0\}$ and $\{e^{tA_r}:t\geq 0\}$ 
on a Banach space 
$\bB$, such that the intersection of the domains 
of the generators 
is dense in $\bB$.  
Clearly, in the context of  
Cauchy problems Trotter's formula states 
the convergence of the
splitting-up approximations 
defined by the splitting 
$$
\frac{\partial v_r}{\partial t}=A_rv_r(t), \qquad r=1,2,...,d_1
$$
to the solution 
of the abstract Cauchy problem 
$$
\frac{\partial u}{\partial t}=Au(t), 
\quad t\geq0,
\quad
u(0)=z.
$$

Our main interest in the present paper is to
increase the accuracy of the splitting-up 
approximations for equation (\ref{17.03.09.03}). 
It is known that the error of the 
splitting-up approximations is 
proportional to $\delta$, the step -size. 
There are, 
however, modifications of these 
approximations which are more
accurate. A celebrated example is 
the Strang symmetric scheme  
$$
u_n(t^{(n)}_i):=(\bS_{\delta/2}^{(1)}
\bS_{\delta/2}^{(2)}...
\bS_{\delta/2}^{(d_1)}
\bS_{\delta/2}^{(d_1)}...
\bS_{\delta/2}^{(2)}\bS^{(1)}_{\delta/2})^iu_0, 
\qquad i=0,1,...,n,
$$
whose error is proportional to $\delta^2$. 
This approximation scheme is presented 
in \cite{S1}, \cite{S3}. 
Other symmetric schemes and their 
generalizations,  
are given in \cite{S1}, \cite{S2}, and \cite{G}. All
these  schemes are of second order accuracy.  
Inspired by 
the above example, for given $k\geq2$ 
one looks for a
composition of splittings    
\begin{equation}                                      \label{15.07.09.03}
\prod_{i=1}^{m}\prod_{j=1}^{d_1}
\bS_{c^{ij}\delta}^{(j)}, 
\end{equation} 
with real numbers $c^{ij}$ and integer 
$m\geq 1$ to be determined, 
such that 
$$
u(\delta)-\prod_{i=1}^{m}\prod_{j=1}^{d_1}
\bS_{c^{ij}\delta}^{(j)}u_0, 
$$
the local error of the corresponding approximation is 
proportional to $\delta^{k+1}$ in appropriate norms. 
Such local error leads to a global error, proportional 
to $\delta^k$, i.e., composition (\ref{15.07.09.03}), 
represents a method of (at least) {\em order\/} $k$.   
The conditions on the numbers $c^{ij}$ and $m$ which  
lead to splitting methods of 
high order have been studied intensively 
in the literature. 
Such methods are obtained in 
\cite{N} for Hamiltonian systems by 
the Baker-Campbell-Hausdorff 
formula. Variations of the
Trotter formula and  the Baker-Campbell-Hausdorff formula 
are used for linear and for nonlinear equations, respectively, 
to show the existence of methods of any order (see 
\cite{MQ},
\cite{Su1},   \cite{TS}, \cite{Y} 
and the literature therein). 
An adaptation of the method of {\em rooted trees} from the 
theory of Runge-Kutte approximations is used in 
\cite{MS}. By \cite{Sh} and 
\cite{V}, however, the numbers
$c^{ij}$   in each scheme (\ref{15.07.09.03}) of order $k\geq 3$ cannot 
be all non-negative. 
Thus, by \cite{Sh} and \cite{V} the above splitting methods
of order greater than or equal to 3 cannot be used to 
approximate the solution of partial differential
equations of parabolic type. 
As R.I.~McLachlan and G.R.W.~Quispel write
on page 392 of \cite{MQ}: ``...splitting was proposed as a cheap
way to retain unconditional stability. Methods with backward time steps
can only be conditionally stable; this stumbling block held up the
development of high-order compositions for years." 

Then the natural question arises, 
as to whether there exists, in the case of parabolic 
equations, a different way from the multiplicative one 
to accelerate the convergence to a higher order. 
 One of our main results consists of
 showing that using the   step size of order $\delta$, 
but organizing the
computations differently, it is indeed
possible to achieve the  accuracy of order
$\delta^{k}$ for any $k$, even if $A_{r}$ are
(degenerate) elliptic operators with coefficients depending on time.
In a subsequent article we intend to show that
our method is much  more universal 
in the sense that it covers very many situations in which
method \eqref{15.07.09.03} works and requires
approximately the same amount of work.

In the present paper we use linear combinations of 
splittings of type (\ref{12.05.09.03}) 
with different step-sizes, 
to achieve arbitrary high accuracy. 
We prove that for any 
given $k\geq0$ there exist 
 absolute constants 
$ b_0,b_1,...,b_k $ expressed by simple formulas 
such that the accuracy of the
approximation 
\begin{equation}                              
                                                 \label{20.09.09.03}
v_n:=b_0u_n+b_1u_{2n}+b_2u_{4n}+...+b_ku_{2^kn}
\end{equation}
is of order $\delta^{k+1}$  
(see Theorem \ref{theorem12.09.09.03} below). 
Here $u_{2^jn}$ is the splitting-up approximation
(\ref{12.05.09.03}) along the grid 
(\ref{grid}), with $2^jn$ in
place of $n$. 
In particular, if $k=1$, we have to deal with
two step sizes: $\delta$ and $\delta/2$, and we get the order
of accuracy $\delta^{2}$. The Strang formula giving the same
order of accuracy, generally,
also requires working with step size $\delta/2$. By the way,
if $A=A_{1}+A_{2}$ and we construct our splitting-up
scheme according to $A=(1/2)A_{1}+A_{2}+(1/2)A_{1}$,
then our approximations just coincide  with the Strang one
and there is no need to use linear combinations to get
the error of order $\delta^{2}$.
 It is also worth noting that the above
coefficients ${\boldsymbol b}=(b_{0},...,b_{k})$
 are given
by 
$\mathbf b:={\mathbf e}_1 V^{-1}$, where 
${\mathbf e}_1:=(1,0,0,...,0)$ and  
$V^{-1}$ is the inverse of
the $k+1\times k+1$ Vandermonde matrix  
$V^{ij}=2^{-(i-1)(j-1)}$.

Our work is of purely theoretical nature and, as the referees
pointed out, much work yet needs to be done before
our results could be used in practical applications.
We restricted ourselves to making the first step
in attacking Problem 10 on page 492 of
\cite{MQ}: ``For systems that evolve
in a semigroup, such as the heat equation,
develop effective methods of order higher than 2''.
However, our results show that each time when one has
{\em any\/} algorithm of implementing standard splitting-up method
to approximating the solutions of the Cauchy problem for
degenerate parabolic equations with sufficiently
smooth coefficients and free terms, one can
improve the rate of convergence to any degree.
For instance, we believe that
usually in practice
one is not doing computations with only one
step size, and  we show  that having,
say, three different step sizes each of which
is of order $\delta$  of accuracy, and just
taking a linear combination of the results, 
one gets an approximation with error
of order $\delta^{3}$.

We have to admit that we do not know if
our methods can be carried over to quasilinear equations
or to equations in domains. In this connection we note that
there is a very active  area of
developing and applying in practice splitting-up methods for
degenerate nonlinear convection-diffusion equations
(see, for instance, \cite{EK} and 185 references therin).
Our equations can be viewed as  belonging to this
area only if $a^{ij}$ are constant. 
However, it is perhaps worth mentioning that
our methods can be applied to solving systems
of (nonlinear) ODEs and we are in the process
of working on this subject.

 Inspired by Richardson's method 
we obtain our results by expanding 
the error 
$u-u_n$ of the splitting-up approximation (\ref{12.05.09.03}) in
powers of 
$\delta=T/n$. This is Theorem \ref{theorem 11.26.06}, the main
theorem of our paper. We use this expansion with   
$\delta=T/2^{j}n$, $j=0,1,2,...,k$, and choose the 
above coefficients 
$b_0,b_1,...,b_k$ to eliminate the terms of order less than $k+1$ in the 
linear combination (\ref{20.09.09.03}). 

The main theorem of the present paper is proved by exploiting  
a new approach of \cite{KG1} and \cite{KG2} to splitting-up 
methods. As we discussed above, the  splitting-up approximation 
(\ref{12.05.09.03}) means that to get the approximation 
at $t+\delta$ from that at $t$, one goes 
back and forth in time $d_1$-times while solving 
equations (\ref{19.03.09.03}), $r=1,2,...,d_1$, successively. 
A basic idea of \cite{KG1} is 
to arrange the splitting continuously in forward time direction,  
and to synchronize with it the original 
equation by time scaling. In this way 
we have differential equations for 
the rearranged splitting-up  
approximations and for the 
time-scaled solution of the original 
equation, 
which enables us
to use methods of the theory of partial differential equations
and not semigroup theory and   get an
expansion  for their difference in terms of powers of $\delta$
even if the coefficients depend on time.
The method  of \cite{KG1} and \cite{KG2} appeared
in connection with splitting-up for stochastic
partial differential equations. It is worth mentioning that
most likely it is impossible to accelerate 
the splitting-up method in this more
complicated situation.

The paper is organized as follows.
In the next section we introduce our
general setting but state the results, 
Theorems \ref{theorem12.09.09.03} 
and \ref{theorem 11.26.06} 
only for the case of time independent
data for the sake of simplicity of presentation.  
Theorem \ref{theorem 11.26.06} 
is proved immediately after its formulation 
on the basis of Theorem
\ref{theorem12.09.09.03}, 
which in turn is proved 
in Sec.~\ref{section 7.27.1}, 
after we prepare some auxiliary facts in
Sec.~ \ref{section 17.01.09.03}. 
In Sec.~ \ref{section 18.01.09.03}
we generalize Theorems \ref{theorem12.09.09.03} 
and \ref{theorem 11.26.06} 
for time dependent data 
and derive some consequences valid in the time-homogeneous
case as well.

In conclusion we introduce some notation used everywhere
below. Throughout the paper
$d\geq1, d_{1}\geq2$ are fixed positive integers,
$K,T$ are fixed finite positive constants, and
$$
D_{i}:=\partial/\partial x^{i},\quad
D_{ij}:=\partial^{2}/\partial x^{i}\partial x^{j},
\quad D_{t}:=\partial/\partial t.
$$
We denote by $W^m_p$ the Sobolev space defined 
as the closure of $C_0^{\infty}$ functions 
$\varphi:\bR^d\to\bR$ in the norm 
$$
\|\varphi\|_{m,p}
:=\Big(\sum_{|\gamma|\leq m}
\int_{\bR^d}|D^{\gamma}\varphi(x)|^p\,dx\Big)^{1/p}, 
$$
where $D^{\gamma}:=D_1^{\gamma_1}...D^{\gamma_d}_d$ 
for multi-indices $\gamma=(\gamma_1,...,\gamma_d)$ 
of length $|\gamma|:=\gamma_1+\gamma_2+...+\gamma_d$. 
Unless otherwise indicated, we use 
the summation convention 
with respect to
repeated indices.  
     
We express our
sincere gratitude to the referees for useful
criticism which helped improve the presentation.

\mysection{Formulation of the main results.
The case of time independent
coefficients}

We consider the problem
\begin{equation}                                  \label{12.09.01}
 D_{t}  u(t,x)
=Lu(t,x)+f(t,x),
\qquad t\in(0,T],\,\,\,x\in\bR^d,
\end{equation}
\begin{equation}                                  \label{13.09.01}
u(0,x)=u_0(x), \quad
x\in\bR^d,
\end{equation}
where $L$ is an operator of the form
$$
L =a^{ij}(t,x)D_{ij}+a^{i}(t,x)D_{i}+a(t,x),
$$
$f$ and $u_0$ are real functions of
$(t,x)\in (0,T]\times\bR^d$ and of
$x\in\bR^d$, respectively. 
We assume that the coefficients
$a^{ij}, a^i, a$ and the derivatives 
$a^{ij}_{x^k}$ of
$a^{ij}$ are bounded Borel functions of 
$(t,x)$. We fix
$p\geq 2$ and assume
that $u_{0}$ and $f$ are measurable and 
$|u_0|^p$ and $|f|^p$ are
integrable  over $\bR^d$ and
over $[0,T]\times\bR^{d}$, respectively.
\begin{definition}
                                            \label{definition 7.28.1}
  By a solution of problem
(\ref{12.09.01})-(\ref{13.09.01})
we mean an
$W^1_p$-valued weakly continuous 
function $u(t)=u(t,\cdot)$ defined
on $[0,T]$ such that for all $\phi\in C_0(\bR^d)$
and $t\in[0,T]$ 
$$
(u(t,\cdot),\phi)=(u(0,\cdot),\phi) +\int_{0}^{t}[-
(a^{ij}D_{i}u(s),D_{j}\phi)
$$
$$
+
((a^{i}-a^{ij}_{x^{j}})D_{i}u(s)+au(s)
+f_{r}(s),\phi)]\,ds,
$$
  where $(\,,\,)$ denotes the
usual inner product in $L^2(\bR^d)$.
Quite often we write equation (\ref{12.09.01})
and similar equations in the form
$$
du(t)=(Lu(t)+f(t))\,dt
$$
bearing in mind the differential of  $u$ in $t$ only. 
\end{definition}
Suppose that we split 
equation (\ref{12.09.01}) 
into the equations
\begin{equation}                                          \label{9.15.11}
D_{t} v(t,x)
=L_rv(t,x)+f_r(t,x),
\qquad t\in(0,T],\,\,\,x\in\bR^d
\end{equation}
with
$$
L_r :=a^{ij}_r(t,x)D_{ij}+a^{i}_r(t,x)D_{i}+a_r(t,x),
\quad 
L=\sum_{ r=1 }^{ d_1 }L_r,
\quad  f=\sum_{ r=1 }^{ d_1 }f_r,
$$
such that these equations are more pleasant 
from point of view of
numerical methods than the original one.
This motivates the multi-stage splitting
method, which we describe below. 
First we need some
assumptions.

Fix an integer 
$l\geq1$. 

\begin{assumption}[ellipticity of $L_r$]
                                    \label{assumption 10.09.01}
For each $r=1,2,...,d_1$ for $dt\times dx$-almost
every $(t,x)\in[0,T]\times\bR^d$
$$
a^{ij}_r(t,x)
\lambda^{i}\lambda^{j}\geq0
$$
for all 
$(\lambda^1,\lambda^2,...,\lambda^d)\in\bR^{d}$.
\end{assumption}

\begin{assumption}                           
                                \label{assumption 10.11.09.03} 
(i) The partial 
derivatives  
$$
 D^{s}_{t} D^{\rho}a^{ij}_r, 
\quad 
 D^{s}_{t} D^{\rho}a_r^{i}, 
\quad
 D^{s}_{t} D^{\rho}a_r
\quad {\text{\rm for}}\,\,i,j=1,2,...,d,\, r=1,2,...,d_1 
$$
exist and by magnitude 
are bounded by $K$ 
for all integers $s\geq 0$ and multi-indices $\rho$, 
satisfying $2s+|\rho|\leq l$.  

(ii) For every integer $s\in[0,l/2]$    
$$
\sup_{t\in[0,T]}\|D^{s}_{t} f_r(t)\|_{l-2s,p}\leq K.
$$ 

(iii) We have $
u_{0}\in W^{l}_{p}$ and $\|u_{0}\|_{l,p}\leq K$.

\end{assumption}

It is well-known that under the above 
conditions equations (\ref{12.09.01}) 
and (\ref{9.15.11}) with initial condition 
$u(0)=u_0$ admit a unique generalized solution 
$u$ and $v$, respectively, 
which are $W^l_p$-valued weakly continuous functions 
of $t\geq 0$
(see, for instance, Theorem \ref{theorem28.11.01} below).  
We want to approximate the solution
$u$, by using the splitting-up method, 
i.e., by solving equations (\ref{9.15.11}) 
successively 
with appropriate initial conditions 
on appropriate time intervals. 
Let us formulate now our splitting-up scheme 
in the case when
the coefficients $a^{ij}_r,a^i_r,a_r$ 
and free terms $f_r$ are independent
of the time variable $t$.

Set $T_n:=\{t_i:=iT/n:i=0,1,2,...,n\}$,
$\delta:=T/n$ for an integer $n\geq 1$. 
Then for fixed $n$
we approximate the
solution $u$
of (\ref{12.09.01})-(\ref{13.09.01}) at 
$t_i=iT/n$ recursively by
$u_{n}(0):=u_0$, 
\begin{equation}                                    \label{switching}
u_{n}(t_{i+1}):=\bS_{\delta}^{(d_1)}
...\bS_{\delta}^{(2)}
\bS_{\delta}^{(1)}
u_{n}(t_i),
\quad i=0,1,2,...,n-1,
\end{equation}
where
$\bS_{t}^{(r)}\psi:=v(t)$
denotes the solution
of equation (\ref{9.15.11}) 
for $t\geq 0$ with initial condition
$v(0)=\psi$. 

It is known that if 
Assumptions \ref{assumption 10.09.01}, 
\ref{assumption 10.11.09.03}  
are satisfied with $ l= m+4$, then 
$$
\max_{t\in T_n}\|u(t)-u_n(t)\|_{m,p}\leq N/n
$$
for all $n\geq 1$, where $N$ depends 
only on $d,d_1,T,K,p,m$. 
Moreover, 
this rate of convergence is sharp 
(see \cite{KG2}, where this result is 
a special case of the rate of 
convergence estimates 
for stochastic PDEs). 
In the present paper we want to show 
that by suitable combinations 
of splitting-up approximations we 
can achieve as fast convergence as 
we wish. We show this by the aid of 
the following theorem on expansion 
of $u_n$ in powers of the step-size 
$\delta$.

\begin{theorem}
                                         \label{theorem 11.26.06}
Let $m\geq0$ and $k\geq 0$ 
be integers. Let 
Assumptions \ref{assumption 10.09.01} and
\ref{assumption 10.11.09.03}  
hold with 
\begin{equation}
                                        \label{10.4.5}
 l\geq 4+m+4k.
\end{equation}
Suppose that 
the coefficients $a^{ij}_r,a^i_r,a_r$ 
and the free terms $f_r$ do not depend on $t$. 
Then  
 for all $n\geq 1$ and  $t\in T_n$ and
$x\in\bR^d$, the following representation holds
$$
u_n(t,x)=u(t,x)
+\delta u^{(1)}(t,x)
$$
\begin{equation}
+\delta^2u^{(2)}(t,x)+...
+\delta^k u^{(k)}(t,x)+R^{(k)}_n(t,x) ,      \label{03.01.02}
\end{equation}
where the functions 
$u^{(1)}$,...,$u^{(k)}$, and 
$R^{(k)}_n$, defined on $[0,T]$,
 are $W^m_p$-valued and
weakly continuous. Furthermore,
$u^{(j)}$, $j=1,2,...,k$, 
are independent of $n$, and
\begin{equation}                            \label{remainder}
\sup_{t\in T_{n}}
\|R^{(k)}_n(t)\|_{m,p}\leq N\delta^{k+1}
\end{equation}
for all $n$, where $N$ depends only on $l$,
$k,d,d_1,K,m,p,T$.
\end{theorem}

\begin{remark}
If $k=0$ and $p=2$, the result holds under
a weaker restriction on $l$: $l\geq3+m$
(see, for instance \cite{KG1}).
For general $p\geq2$ and $k=0$ the result
is proved  in \cite{KG2}.
\end{remark}

We prove Theorem \ref{theorem 11.26.06}
 in Section \ref{section 7.27.1}. 
Now we deduce from it
a result on the acceleration of 
the splitting-up method. 
Let $V$ denote the square matrix defined by 
$V^{ij}:=2^{-(i-1)(j-1)}$, $i,j=1,...,k+1$.  
Notice that the determinant of $V$ 
is the Vandermonde
determinant, generated by 
$1,2^{-1},...,2^{-k}$, 
and hence it is
different from 0. 
Thus $V$ is invertible. Set 
$\mathbf b:=(b_0,b_1,...,b_k)
:=(1,0,0,...,0)V^{-1}$, and define 
$$
v_n(t):=\sum_{j=0}^{k}b_ju_{2^jn}(t),
\quad t\in T_n:=\{iT/n: i=0,1,..,n\}, 
$$
where $u_{2^jn}$ is the splitting-up 
approximation
based on the grid 
$T_{2^jn}:=\{iT/(2^jn):i=0,1,..,2^jn\}$.

\begin{theorem}                                  \label{theorem12.09.09.03}
Let $m\geq 0$ and $k\geq 0$ be any integers. 
Let 
Assumptions  
\ref{assumption 10.09.01} and  
\ref{assumption 10.11.09.03}
hold with    $l$ satisfying (\ref{10.4.5}).
Suppose that the coefficients  
$a^{ij}_r,a^i_r,a_r$ 
and the free terms $f_r$ do not depend on $t$. 
Then 
$$
 \max_{t\in T_n}\|v_n(t)-u(t)\|_{m,p}
\leq N\delta^{k+1}, 
$$
where $N$ is a constant, depending only on $l$,
$k,d,d_1,K,m,p,T$. 
\end{theorem}
\begin{proof}
By Theorem \ref{theorem 11.26.06} 
$$
u_{2^jn}=
u+\sum_{i=1}^k\frac{\delta^i}{2^{ji}}
u^{(i)}+R^{(k)}_{2^jn},
\quad j=0,1,...,k.
$$
Therefore for all $n\geq 1$
$$
v_n=\sum_{j=0}^kb_ju_{2^jn}
=(\sum_{j=0}^{k}b_j)u
+\sum_{j=0}^{k}\sum_{i=1}^kb_j
\frac{\delta^i}{2^{ij}}u^{(i)}
+\sum_{j=0}^kb_jR^{(k)}_{2^{j}n}
$$
$$
=u+\sum_{i=1}^k\delta^iu^{(i)}
\sum_{j=0}^{k}\frac{b_j}{2^{ij}}
+\sum_{j=0}^kb_jR^{(k)}_{2^{jn}}
=u+\sum_{j=0}^kb_jR^{(k)}_{2^{j}n},
$$
since $\sum_{j=0}^{k}b_j=1$ and  
$\sum_{j=0}^{k} b_j 2^{-ij} =0$ for   
 $i=1,2,...k$ by the definition of $(b_0,...,b_k)$. 
Hence $v_n-u=\sum_{j=0}^kb_jR^{(k)}_{2^{j}n}$, and 
$$
\max_{t\in T_n}\|v_n(t)-u(t)\|_{m,p}
=\max_{t\in T_n}
\|\sum_{j=0}^kb_jR^{(k)}_{2^{j}n}(t)\|_{m,p}
$$
$$
\leq \sum_{j=0}^k|b_j|
\max_{t\in T_n}\|R^{(k)}_{2^{j}n}(t)\|_{m,p}
\leq N\delta^{k+1}, 
$$
by (\ref{remainder}), where $N$ is a constant 
depending only on $l$, $T$, $K,d,d_1,m,p,k$. 
\end{proof}

\begin{remark}Assume that $u^{(1)}=0$ in expansion 
(\ref{03.01.02}). This happens, for example, for Strang's 
splitting, which is a special case of our splitting-up 
scheme, as it is explained in the Introduction. In this case 
we need only take $k$ terms in the linear combination to achieve accuracy 
of order $k+1$. Namely, we define now $v_n(t)$ by 
$$
v_n(t):=\sum_{j=0}^{k-1}\lambda_ju_{2^jn}(t), \quad t\in T_n, 
$$
where 
$$
(\lambda_0,\lambda_1,...,\lambda_{k-1})
:=(1,0,...,0)V^{-1}, 
$$
and $V$ is now a $k\times k$ Vandermond matrix with 
entries $V_{i1}:=1$, 
$V_{i,j}:=2^{-(i-1)j}$ for $i=1,2,...,k$ and  
$j=2,...,k$. Then Theorem \ref{theorem12.09.09.03} remains 
valid, what one can prove in the same way as Theorem 
\ref{theorem12.09.09.03} is proved. 
For example, 
$$
v_n(t):=-\frac{1}{3}u_n(t)+\frac{4}{3}u_{2n}(t), \quad t\in T_n
$$
is an approximation of accuracy $\delta^3$ in
{\color{blue} the} case of Strang's 
splitting. 
\end{remark}

\mysection{Auxiliary Results}         \label{section 17.01.09.03}    
Let us consider the partial 
differential equation
\begin{equation}
                                     \label{14.03.07.02}
du(t,x)=(Lu(t,x)+f(t,x))\,dA(t),
\qquad t\in(0,T],\,\,\,x\in\bR^d,
\end{equation}
\begin{equation}                     \label{15.03.07.02}
u(0,x)=u_0(x), \quad
x\in\bR^d,
\end{equation}
where $L$ is an operator of the form
$$
L =a^{ij}(t,x)D_{ij}+a^{i}(t,x)D_{i}+a(t,x),
$$
$A=A(t)$ is a continuous increasing function
starting from $0$,
$f$ and $u_0$ are real functions of
$(t,x)\in (0,T]\times\bR^d$ and of
$x\in\bR^d$, respectively.
Fix  an integer $l\geq0$ and a real number $p\geq2$. We 
understand the solution
in the spirit of Definition \ref{definition 7.28.1} and 
make the following assumptions.

\begin{assumption}[smoothness of the coefficients]
                                 \label{assumption 16.04.07.02}
The coefficients of $L$
are measurable. The derivatives 
in $x\in\bR^{d}$ of the 
coefficients $a^{ij}$ 
up to order $2\vee l$, 
  of the coefficients 
$a^{i}(t,x)$ up to 
order
$1\vee l$,   and of 
$a(t,x)$  up to order
$l$ exist
for any $t\in(0,\infty)$, 
and by magnitude are bounded by $K$.
\end{assumption}
\begin{assumption}
                                    \label{assumption 06.01.02}
We have
$$
u_{0}\in W^{l}_{p},\quad f \in 
L_{p}([0,T],W^{l}_{p}).
$$
\end{assumption}

\begin{assumption}[ellipticity of $L$]
                                    \label{assumption 17.04.07.02}
For all $t\geq0$,
$x\in\bR^{d}$, and
$\lambda\in\bR^{d}$,
we have 
$$
a^{ij}(t,x)
\lambda^{i}\lambda^{j}\geq0.
$$
\end{assumption}

\begin{assumption}
                                   \label{assumption 03.07.02}
The function $A$ is absolutely continuous and
$$
\dot A(t):=\frac{d}{dt}A(t)\leq K 
$$
for $dt$-almost every $t\geq0$. 
\end{assumption}

The following result is well-known in PDE theory    
(after replacing $dA $ in (\ref{14.03.07.02}) with 
$\dot A \,dt$ we easily get it, for instance, from \cite{KR} or from
  Theorem 3.1 in 
\cite{KG1}).  

\begin{theorem}                                 \label{theorem28.11.01} 
Under Assumptions 
\ref{assumption 16.04.07.02},
\ref{assumption 06.01.02},
\ref{assumption 17.04.07.02},
and
\ref{assumption 03.07.02} with 
 $l\geq1$  
  the Cauchy problem
(\ref{14.03.07.02})-(\ref{15.03.07.02}) has a unique
generalized solution $u$. 
  If 
Assumptions 
\ref{assumption 16.04.07.02},
\ref{assumption 06.01.02},
\ref{assumption 17.04.07.02},
and
\ref{assumption 03.07.02} 
hold with $l\geq0$, and $u$ is a generalized  solution 
of (\ref{14.03.07.02})-(\ref{15.03.07.02}), then 
for every integer $l_1\in[0,l]$  
$$
\sup_{t\in[0,T]}\|u(t)\|^{p}_{l_1,p}\leq
N\Big\{\|u_0\|^{p}_{l_1,p}
+\int_0^T\|f(t)\|^{p}_{l_1,p}\,dt\Big\},
$$ 
where $N$ is a constant depending only on $T$,
$K,l,p,d$. 
\end{theorem}

Under the assumptions of 
Theorem \ref{theorem28.11.01} 
  let $\cR f$ denote the solution of equation 
(\ref{14.03.07.02}) with
initial data $u_0=0$. Then by virtue of 
Theorem~\ref{theorem28.11.01}
$$
\cR:L_p([0,T],W_p^l)
\to C_w([0,T],W_p^l)
$$
is a bounded linear operator,
where $C_w([0,T],W_p^l)$
denotes the Banach space of 
weakly continuous $W_p^l$-valued functions
$u=u(t)$, $t\in[0,T]$ with the
norm
$
\sup_{t\in[0,T]}\|u(t)\|_{l,p}.
$

\smallskip
Let us now consider the equation
\begin{equation}
                                                     \label{16.27.08}
du(t,x)=Lu(t,x)\,dA(t)+g(t,x)\,dH(t),
\end{equation}
$$
(t,x)\in(0,T]\times\bR^d,
$$
where
$g$ is a real-valued function of
$(t,x)\in [0,T]\times\bR^d$
and $H$ is an
absolutely continuous function of $t\in[0,T]$.
\begin{assumption}                            \label{assumption 10.04.07.02}
We have 
$
g\in L_p([0,T],W^{l+2}_p)
$, 
and there exists 
$
g^{\prime}\in L_p([0,T],W^{l}_p)
$
such that 
$$
d(g(t),\phi)=(g^{\prime}(t),\phi)\,dA(t),
\quad t\in[0,T]
$$
for all $\phi\in C_0(\bR^d)$.
\end{assumption}

\begin{lemma}                               \label{lemma16.27.08}
Under Assumptions
\ref{assumption 16.04.07.02},
\ref{assumption 17.04.07.02},
\ref{assumption 03.07.02}, 
\ref{assumption 10.04.07.02}  
  with $l\geq 1$  
equation
(\ref{16.27.08})
with zero initial data has
a unique generalized solution $u$.
Moreover,
\begin{equation}                             \label{12.03.06.02}
u=\cR(H(Lg-g^{\prime}))+Hg=:\cQ(H,g). 
\end{equation}
 If Assumptions
\ref{assumption 16.04.07.02},
\ref{assumption 17.04.07.02},
\ref{assumption 03.07.02}, 
\ref{assumption 10.04.07.02}  
hold with $l\geq0$ and 
equation (\ref{16.27.08})
with zero initial data admits 
a generalized solution $u$, then 
for every integer $l_1\in[0,l]$ 
\begin{multline}                              \label{18.27.08}
\sup_{t\in[0,T]}
\|u(t)\|_{l_1,p}
\leq
N\sup_{t\in[0,T]}|H(t)|
\Big(\sup_{t\in[0,T]}\|g(t)\|_{l_1,p}\\
+\{\int_0^T(\|g(t)\|^p_{l_1+2,p}
+\|g^{\prime}(t)\|^p_{l_1,p})\,dt\}^{1/p}\Big),
\end{multline}
where
$N$ is a
constant depending only on
$p,d,K,l,T$.
\end{lemma}
\begin{proof}
Note that 
$$
d((g(t),\phi)H(t))=(g(t),\phi)\,dH(t)
+(g^{\prime}(t),\phi)H(t)\,dA(t),
$$
for all $\phi\in C_0(\bR^d)$. Therefore 
$u$ solves equation  (\ref{16.27.08}), when
$w(t,x):=u(t,x)-H(t)g(t,x)$ solves
$$
dw(t,x):=\{Lw(t,x)
+H(t)(Lg(t,x)-g^{\prime}(t,x))\}\,dA(t).
$$
Hence equality (\ref{12.03.06.02}) follows by
Theorem \ref{theorem28.11.01},
and it implies (\ref{18.27.08}).
\end{proof}

\begin{remark}
                                                  \label{remark 10.12.1}
We often consider equation 
(\ref{16.27.08}) when  
$d H(t)/dt$ 
is bounded. 
Then under Assumptions 
\ref{assumption 16.04.07.02},
\ref{assumption 17.04.07.02},
\ref{assumption 03.07.02}, 
\ref{assumption 10.04.07.02}  
with $l\geq1$,   
equation (\ref{16.27.08})
with zero initial data has a unique 
generalized solution $u$ by 
Theorem \ref{theorem28.11.01}.
By Theorem \ref{theorem28.11.01}
this solution belongs to $C_{w}([0,T],W^{l }_{p})$
and its norm in this space admits an estimate
with a constant depending on the bound for $\dot H$
and only the $L_{p}([0,T],W^{l}_{p})$-norm of $g$.
It is important that the function $H$ enters
(\ref{18.27.08}) only through $\sup|H|$ and not any
characteristic of its derivative,
however for that we pay a price
requiring $g$ to have more derivatives.
\end{remark}
 
\mysection{Proof of Theorem \ref{theorem 11.26.06}}
                                  \label{section 7.27.1}

Throughout this section 
the assumptions of Theorem \ref{theorem 11.26.06} are supposed to be
satisfied. In particular, $l\geq4$.
Fix $n$ and introduce $\delta=T/n$. 
We use the idea from \cite{KG1} and \cite{KG2} 
of rearranging the splitting method in forward time. 
We  achieve this by considering the equation 

\begin{equation}                \label{equation 09.26.06}
dw(t,x)=
\sum_{r=1}^{d_1}(L_rw(t,x)+f_r)\,dA_r(t),\quad
w(0,x)=u_0(x), 
\end{equation}
where the time change $A_r$,
$r=1,...,d_{1}$, is defined by 
the requirements that $A_{r}(0)=0$,
$A_{r}(t)$ be absolutely continuous,
and its derivative in time $\dot A_r$ be periodic
with period $d_{1}\delta$ and
\begin{equation}
                                                   \label{9.24.1}
\dot A_r(t)
=1_{[r-1,r]}(t/\delta), 
\quad t\in[0,d_1\delta]\quad(\mbox{a.e.)}.
\end{equation}

 Instead of the original 
Cauchy problem (\ref{12.09.01})-(\ref{13.09.01}) 
we consider 
$$
dv(t,x)=
(Lv(t,x)+f)\,d A_{0}(t),\quad
v(0,x)=u_0(x),
$$
where 
\begin{equation}                                \label{10.17.09.03}
A_{0}(t):=t/d_1 .
\end{equation} 
Clearly, $v(t)=u(A_{0}(t))$, and 
$$
v(d_1t)=u(t), \quad w(d_1t)=u_{n}(t) 
\quad\text{for all }t\in T_n. 
$$
Therefore our aim is to show
that Theorem \ref{theorem 11.26.06} holds 
with $v$ and $w$ in place of $u$, 
and $u_{n}$, respectively, 
for all $t=id_1\delta$, $i=0,1,...n$.  
To this end first we introduce some notation.  
We call a sequence of numbers 
$\alpha=\alpha_1\alpha_2...\alpha_{i}$ 
a multi-number of length $|\alpha|:=i$, 
if 
$\alpha_j\in\{0,1,2,...,d_1\}$. 
The reader should notice the difference
between multi-numbers and multi-indices. 
The set of all multi-numbers is denoted
by $\cN$. 
For every multi-number 
$
\alpha 
$ 
we
define a function
$B_{\alpha}:[0,\infty)\to\bR$ and a 
number $c_{\alpha}$ recursively starting 
as follows: 

\begin{equation}                                    \label{21.16.09.03} 
\quad B_{\gamma}:=\delta^{-1}(A_{\gamma}-A_{0}),  
\quad  c_{\gamma}=0 
\qquad \text{for }\gamma=0,1,2,..d_1.   
\end{equation} 
If for every multi-number
$\beta=\beta_1...\beta_i$ of length $i$ 
the function 
$B_{\beta}$ and the number 
$c_{\beta}$ are defined, then
\begin{equation}
                                                  \label{18.19.09.03}
c_{\beta\gamma}:={\delta}^{-1}
\int_0^{d_1\delta}B_{\beta}(s)
\dot A_{\gamma}(s)\,ds,
\end{equation}

\begin{equation} 
                                                  \label{22.16.09.03}
B_{\beta \gamma}(t)
:=\delta^{-1}
\int_0^t(B_{\beta}(s)\dot A_{\gamma}(s)
-c_{\beta\gamma}\dot A_{0}(s))\,ds
\end{equation}
for 
$\gamma=0,1,2,...,d_1$, where 
$\dot A_{\gamma}(s):=d A_{\gamma}(s)/ds$. 

Notice that by (\ref{22.16.09.03}) we have 

\begin{equation}                                       \label{16.18.09.03}
B_{\beta}(t)\,dA_{\gamma}(t)=
c_{\beta\gamma}\,dA_{0}(t)+\delta\,dB_{\beta\gamma}(t) 
\end{equation} 
for all multi-numbers 
$\beta$ and 
$\gamma=0,1,2,...,d_1$. 
We will often make use of this equality 
and of the following lemma. 

\begin{lemma}                                 \label{lemma10.25.06}
For every $\alpha\in\cN$ 
the function $B_{\alpha}$
is $d_1\delta$-periodic, i.e.,
$B_{\alpha}(t+d_1\delta)=B_{\alpha}(t)$ 
for all $t\geq0$,  and
$B_{\alpha}(id_1\delta)=0$ for
all integer
$i\geq0$. Moreover,
the numbers $c_{\alpha}$, the functions 
$C_{\alpha}(t):=B_{\alpha}(\delta t)$,
 and
$$
 \sup_{t\geq0}|B_{\alpha}(t)|
=\sup_{t\geq0}|C_{\alpha}(t)|
$$ 
are finite and do not depend
on $\delta$.
\end{lemma}
\begin{proof}That the first assertion
is true for $\alpha=0,...,d_{1}$ is almost obvious.
If it is true for  $\alpha=\beta$, where $\beta$
is a multi-number,
 then the integrand in (\ref{22.16.09.03})
is $d_{1}\delta$-periodic and  by definition of $c_{\beta\gamma}$
its integral over the period is zero. It follows that
the first assertion holds for $\alpha=\beta\gamma$,
so the induction on
the length $|\alpha|$ finishes the proof
of the first assertion.

To prove the second one we again use the induction on
the length $i=|\alpha|$. This statement is true when
$|\alpha|=1$.
Assume that it is true for all
multi-numbers $\beta$ of
length $i$ and notice that according to
(\ref{9.24.1}) and (\ref{10.17.09.03})
$\dot A_\gamma(\delta s)$ are $d_1$-periodic
in $s$  and
independent of $\delta$. 
Therefore,
$$
c_{\beta\gamma} =\frac{1}{\delta}
\int_0^{d_1\delta}
B_{\beta}(s)\dot A_\gamma(s)\,ds,
=\int_0^{d_1}
C_{\beta}(s)\dot A_\gamma(\delta s)\,ds
$$
is independent of $\delta$ by the induction 
hypothesis. Similar argument works for
$C_{\alpha}(t)$.
\end{proof}

We use the notation 
$\cR f$ and $\cQ_{\alpha} g$, $\alpha\in \cN$, for 
the solutions of equations
(\ref{14.03.07.02})  and (\ref{16.27.08}), respectively,
with zero initial condition and $A_{0}$ 
and  $B_{\alpha}$ in place of  $A $ and
$H$,  respectively. Notice that, unlike in the case
of uniformly parabolic operators, $\cR$ and $\cQ_{\alpha}$
do not increase regularity.

The following lemma exhibits our two main technical tools:
centering $B_{\alpha}$ and
integrating by parts with respect to $t$.  
\begin{lemma}                                   \label{lemma 11.09.07.02}

Take some functions
$$
h \in L_{p}([0,T], W^{1}_{p}),\quad
h_{r}\in L_{p}([0,T], W^{1}_{p}),\quad r=1,...,d_{1},
\quad h_{0}=0.
$$
 Let
$u$ be a solution of the ``equation''
$$
du=\sum_{r=1}^{d_{1}}h_{r}\,dA_{r}
$$
$u(0)\in W^{1}_{p}$, 
which is a particular
case of equation (\ref{equation 09.26.06}) when $L_{r}\equiv0$.
Finally, let $Lu\in L_{p}([0,T], W^{1}_{p})$.
Then  for any  
  $\alpha\in\cN$ 
\begin{equation}
                                             \label{9.23.1}
\cR(B_{\alpha} h)=c_{\alpha0}\cR h 
+\delta \cQ_{\alpha0} h, 
\end{equation}
\begin{equation}
                                                     \label{9.23.2}
\cQ_{\alpha}u=\cR(c_{\alpha0}Lu-c_{\alpha r}h_{r})
+\delta\cQ_{\alpha0}Lu -\delta\cQ_{\alpha r} h_{r} 
+B_{\alpha}u.
\end{equation}

\end{lemma}
\begin{proof} 
To prove (\ref{9.23.1}) it suffices to use
the definitions of $\cR$ and $\cQ_{\beta}$
(see Theorem  \ref{theorem28.11.01} and 
Lemma  \ref{lemma16.27.08}) and use that
by virtue of (\ref{16.18.09.03}) 
for $\varphi=\cR(B_{\alpha}h)$ we have 
$$
d\varphi=L\varphi\,dA_{0}+B_{\alpha}  h\,dA_{0}
=L\varphi\,dA_{0}+c_{\alpha0}  h\,d A_{0}
+\delta h \,dB_{\alpha0}.
$$

To prove (\ref{9.23.2}) observe that by definition
 $\theta:=\cQ_{\alpha}u$ 
satisfies
$$
d\theta=L\theta\,d A_{0}+u\,dB_{\alpha},\quad\theta(0)=0.
$$
This and (\ref{16.18.09.03}) imply  
 that $\psi:=\theta-uB_{\alpha}$ satisfies $\psi(0)=0$ and
$$
d\psi=L\psi\,dA_{0}+LuB_{\alpha}\,dA_{0}
-h_{r}B_{\alpha}\,dA_r 
$$
$$
=L\psi\,dA_{0}+
c_{\alpha0}Lu\,dA_{0}+\delta Lu\,dB_{\alpha0}
 -c_{\alpha r}h_{r}\,dA_{0}-\delta h_{r}
\,dB_{\alpha r}.
$$
Now (\ref{9.23.2}) follows from the definitions
of $\cR$ and $\cQ_{\beta}$. 
The proof of the lemma is complete. 
\end{proof}

Now we 
introduce  some differential 
operators $L_{\gamma}$ and 
functions $f_{\gamma}$ 
defined for 
multi-numbers $\gamma$ as follows: 
$L_0:=0$, $f_0:=0$,  
$$
L_{\gamma}:=L_r,\quad f_{\gamma}:=f_r
$$
for $\gamma=r\in\{1,2,...,d_1\}$, and 
$$
L_{\gamma0}:=LL_{\gamma},
\quad L_{\gamma r}:=-L_{\gamma}L_r
$$
$$
f_{\gamma0}:=Lf_{\gamma},
\quad f_{\gamma r}:=-L_{\gamma}f_{r}
$$
for $r=1,2,..., d_1$.

In this notation we have the following.
\begin{lemma}
                                                  \label{lemma 9.28.1}
Let $\alpha,\beta\in\cN$ and
  $2|\beta|+3\leq l$. Then
\begin{multline}                                
                                                     \label{9.28.1}
\cQ_{\alpha}(L_{\beta}w+f_{\beta})
=c_{\alpha r}\cR(L_{\beta r}w+f_{\beta r})\\
+\delta\cQ_{\alpha r}(L_{\beta r}w+f_{\beta r})+B_{\alpha}
(L_{\beta}w+f_{\beta}).
\end{multline}

\end{lemma}
\begin{proof} 
 It follows from formula (\ref{9.23.2}) applied to
$u:=L_{\beta}w+f_{\beta}$, when 
$h_{r}=L_{\beta}L_{r}w+L_{\beta}f_{r}$
(remember $f_{r}$ are independent of $t$), that the left part of
(\ref{9.28.1}) equals
$$
\cR[c_{\alpha0}(LL_{\beta}w+Lf_{\beta})
-c_{\alpha r}(L_{\beta}L_{r}w+L_{\beta}f_{r})]
$$
$$
+\delta\cQ_{\alpha0}(LL_{\beta}w+Lf_{\beta})-
\delta\cQ_{\alpha r}(L_{\beta}L_{r}w+L_{\beta}f_{r})+B_{\alpha}
(L_{\beta}w+f_{\beta}),
$$
which is easily seen to be equal to the right-hand side
of (\ref{9.28.1}).
 
\end{proof}

We derive from (\ref{9.28.1}) one of the most important
formulas.

\begin{proposition}                           \label{01.24.06}
Let $\kappa\geq 0$ be an integer 
and $l\geq2\kappa+3$.  
 Then 
$$
w=v+\sum_{i=1}^{\kappa}\delta^i
\sum_{|\alpha|=i}B_{\alpha}(L_{\alpha}w
+f_{\alpha})
$$
\begin{equation}                                \label{02.24.06}
+\sum_{i=1}^{\kappa}\delta^i\sum_{|\alpha|=i+1}c_{\alpha }
\cR(L_{\alpha}w+f_{\alpha})
+\delta^{\kappa+1}r ^{(\kappa)}
\end{equation}
for all $t\in[0,d_1T]$, where  
$$
r ^{(\kappa)} =\sum_{|\alpha|=\kappa+1}
\cQ_{\alpha}(L_{\alpha}w+f_{\alpha}). 
$$
\end{proposition}
\begin{proof}
First notice that for $\varphi_0:=w-v$ we have 
$$
d\varphi_0=d(w-v)=L\varphi_0\,d A_{0}+\delta
(L_rw+f_r)\,dB_r, 
$$
which proves (\ref{02.24.06}) for $\kappa=0$. 

Next we fix a $\kappa\geq1$ and transform $r^{(i)}$,
for $i=0,...,\kappa-1$, by
applying (\ref{9.28.1})
 with $\alpha=\beta$ and $|\alpha|=i+1$ when
$f_{r}\in W^{2|\beta|+3}_{p}$.
Then we get
$$
r^{(i)}
=\sum_{|\alpha|=i+1,|\beta|=1}c_{\alpha\beta}\cR(L_{\alpha\beta}w
+f_{\alpha\beta})+\delta\sum_{|\alpha|=i+1,|\beta|=1}
\cQ( L_{\alpha\beta}w
+f_{\alpha\beta})
$$
$$
+\sum_{|\alpha|=i+1}B_{\alpha}(L_{\alpha}w+f_{\alpha})
=\sum_{|\alpha|=i+1}B_{\alpha}(L_{\alpha}w+f_{\alpha})
$$
$$
+\sum_{|\alpha|=i+2}c_{\alpha }\cR(L_{\alpha }w
+f_{\alpha })+\delta r^{(i+1)}.
$$
This shows how $r^{(0)},r^{(1)},...,r^{(\kappa)}$ are related to each other
and certainly proves the proposition.

\end{proof}

Decomposition (\ref{02.24.06}) looks very much like (\ref{03.01.02})
the only difference being that the factors of $\delta^{j}$
depend on the approximating function $w$ and the coefficients
of $\delta^{j}$ in the second term on the right contain
$B_{\alpha}$ which is no  power series in $\delta$.
However, observe that we have to estimate the difference $v-w$
only at the points $id_{1}\delta$ at which all $B_{\alpha}$
vanish.

Our next step is
to ``solve''  (\ref{02.24.06}) with respect to $w$ 
by the method of successive iterations, that is
by   substituting $w$ given by (\ref{02.24.06})
into the right-hand side of the same equation. In the process
of doing so we   encounter only one difficulty
when the second term on the right is plugged into
the third one and we have to develop  expressions
like $\cR(B_{\alpha} u)$ into power series in $\delta$.
We transform these terms by using~(\ref{9.23.1})
and (\ref{9.28.1}).

First we introduce the notation
$$
w_{\beta}=L_{\beta}w+f_{\beta},
$$
observe that in these terms (\ref{9.28.1}) is rewritten as
\begin{equation}
                                                       \label{9.29.1}
\cQ_{\alpha}w_{\beta}=c_{\alpha r}\cR w_{\beta r}
+B_{\alpha}w_{\beta}+\delta\cQ_{\alpha r}w_{\beta r},
\end{equation}
and note the following.
\begin{lemma}
                                                 \label{lemma 9.29.1}
If $\kappa\geq0$ is an integer and $\alpha,\beta\in\cN$
and $2(|\beta|+\kappa)+1\leq l$,
then
$$
\cR(B_{\alpha}w_{\beta})=\sum_{i=0}^{\kappa }
\delta^{i}\sum_{|\gamma|=i}c_{\alpha0\gamma}\cR w_{\beta\gamma}
$$
\begin{equation}
                                                   \label{9.29.2}
+\sum_{i=1}^{\kappa }\delta^{i}\sum_{|\gamma|=i-1}
B_{\alpha0\gamma}w_{\beta\gamma}
+\delta^{\kappa+1}\sum_{|\gamma|=\kappa }\cQ_{\alpha0\gamma}
w_{\beta\gamma},
\end{equation}
where for any multi-numbers $\mu,\nu$
$$
\sum_{|\gamma|=0}c_{\nu\gamma}\cR w_{\mu\gamma}:=
c_{\nu}\cR w_{\mu},\quad
\sum_{|\gamma|=0}B_{\nu\gamma}w_{\mu\gamma}:=
B_{\nu}w_{\mu},
$$
$$
\sum_{|\gamma|=0}\cQ_{\nu\gamma}
w_{\mu\gamma}:=\cQ_{\nu}
w_{\mu}.
$$
\end{lemma}

\begin{proof} If $\kappa=0$,
(\ref{9.29.2}) follows from (\ref{9.23.1}). If $\kappa\geq1$,
by applying repeatedly (\ref{9.29.1}) as in the proof
of Proposition \ref{01.24.06} we find
$$
\cQ_{\alpha}w_{\beta}=\sum_{i=0}^{\kappa-1}
\delta^{i}\sum_{|\gamma|=i+1}c_{\alpha\gamma}\cR w_{\beta\gamma}
$$
$$
+\sum_{i=0}^{\kappa-1}\delta^{i}\sum_{|\gamma|=i}
B_{\alpha\gamma}w_{\beta\gamma}+\delta^{\kappa }
\sum_{|\gamma|=\kappa }\cQ_{\alpha\gamma}w_{\beta\gamma}.
$$
We use this formula for $\alpha0$ in place of $\alpha$
and finish the proof by referring to (\ref{9.23.1}).

\end{proof}

Let $\cM$ denote the set of multi-numbers 
$\gamma_1\gamma_2...\gamma_i$ with 
$\gamma_j\in\{1,2,...,d_1\}$, 
$j=1,2,...,i$, and integers 
$i\geq1$. 
\begin{lemma}                                 \label{lemma 19.20.10.03}
The following 
statements hold.
 
(i) Let $\gamma=\gamma_1\gamma_2...\gamma_i
\in\cM$ be
such that $|\gamma|=i\leq1+l/2$. Then 
$$
L_{\gamma}=(-1)^{|\gamma|-1}L_{\gamma_1}...L_{\gamma_{i}},
\quad f_{\gamma}
=(-1)^{|\gamma|-1}L_{\gamma_1}...
L_{\gamma_{i-1}}f_{\gamma_{i}}, 
$$

(ii) Let $\beta,\gamma\in\cM$ be such that 
$|\beta|+|\gamma|\leq 1+l/2$. Then 
$$
L_{\beta}L_{\gamma}=-L_{\beta\gamma},
\quad
L_{\beta}f_{\gamma}=-f_{\beta\gamma}. 
$$

(iii) Let $\alpha\in\cN$ be
such that $\rho:=|\alpha|\leq1+l/2$. 
Then there 
exist constants
$c(\gamma)=c(\alpha,\gamma)\in\{0,\pm1\}$ 
defined for all $\gamma\in\cM$  with 
$|\gamma|=\rho$, such that 
\begin{equation}                                \label{equation 19.20.10.03}
L_{\alpha}
=\sum_{\gamma\in\cM,|\gamma|=\rho}c(\gamma)L_{\gamma}, 
\quad 
f_{\alpha}
=\sum_{\gamma\in\cM,|\gamma|=\rho}c(\gamma)f_{\gamma}. 
\end{equation}
\end{lemma}

\begin{proof}
Part (i) follows immediately from the definition 
of $L_{\gamma}$, $f_{\gamma}$ by induction on $|\gamma|$. 
Part (i) obviously implies Part (ii). 
Part (iii) clearly holds for $\alpha=0$
and 
$\alpha=r\in\{1,...,d_1\}$. Assume that equations 
(\ref{equation 19.20.10.03}) hold for 
some $\alpha\in\cN$, $|\alpha|<1+l/2$. Then 
$$
L_{\alpha r}
=-L_{\alpha}L_r=-\sum_{|\gamma|=|\alpha|}c(\gamma)L_{\gamma}L_r
=\sum_{|\gamma|=|\alpha|}c(\gamma)L_{\gamma r},
$$
$$
f_{\alpha r}
=-L_{\alpha}f_r=-\sum_{|\gamma|=|\alpha|}c(\gamma)L_{\gamma}f_r
=\sum_{|\gamma|=|\alpha|}c(\gamma)f_{\gamma r}
$$
for $r\in\{1,2...,d_1\}$, and 
$$
L_{\alpha0}=LL_{\alpha}=\sum_{r=1}^{d_1}L_r
\sum_{\gamma\in\cM,|\gamma|=\rho}c(\gamma)L_{\gamma}=
-\sum_{r=1}^{d_1}\sum_{\gamma\in\cM,|\gamma|=\rho}
 c(\gamma)L_{r\gamma}, 
$$
$$
f_{\alpha0}=Lf_{\alpha}=\sum_{r=1}^{d_1}L_r
\sum_{\gamma\in\cM,|\gamma|=\rho}c(\gamma)f_{\gamma}
=-\sum_{r=1}^{d_1}
\sum_{\gamma\in\cM,|\gamma|=\rho} c(\gamma)f_{r\gamma}, 
$$
which prove (iii) by induction on $|\alpha|$. 
\end{proof}

We introduce sequences 
$\sigma=(\beta_1,\beta_2,...,\beta_i)$ of multi-numbers $\beta_j\in\cM$, 
where $i\geq1$ is any integer, and set 
$|\sigma|:=|\beta_1|+|\beta_2|+...+|\beta_i|$. We consider also 
the `empty sequence' $e$ of length $|e|=0$, and denote the set of all 
these sequences by $\cJ$. For $\sigma=(\beta_1,\beta_2,...,\beta_i)$,
$i\geq1$,  
we define 
$$
S_{\sigma}=\cR 
L_{\beta_{1}} \cdot...\cdot\cR L_{\beta_{i}} 
$$ 
 and for $\sigma=e$ we set 
$$
S_e=\cR.
$$
Notice that $S_{\sigma}$ involves $2|\sigma|$
derivatives with respect to $x$ and certain number of operators
$\cR$ which do not increase regularity. Therefore, basically,
$S_{\sigma}$ has the power of a differential operator of
$2|\sigma|$'th order. 
If we have a collection of functions $g_{\nu}$ indexed by a parameter
$\nu$ taking values in a set $A$, then we use the notation
$$
\sumstar_{\nu\in A}g_{\nu}
$$
for any linear combination of $g_{\nu}$ with coefficients
independent of the argument of $g_{\nu}$ and of $\delta$.
For instance,
$$
\sumstar_{  A}S_{\sigma}w_{\gamma}=
\sumstar_{(\sigma,\gamma)\in A}S_{\sigma}w_{\gamma}= 
\sum_{(\sigma,\gamma)\in A}
c(\sigma,\gamma)S_{\sigma}w_{\gamma} ,
$$
where $c( \sigma ,\gamma)$ are certain constants independent of 
$\delta$. These constants are allowed to 
change from one occurrence
to another. 

For functions $u=u(t,x)=u( \delta ,t,x)$ depending on
the parameter $\delta$ we write $u=O_{m}(\delta^{\kappa})$,
if 
$$
\sup_{\delta} \delta^{-\kappa}
\sup_{t\in[0,d_{1}T]}\|u(t)\|_{m,p}<\infty.
$$
 We also use the following sets
$$
A(i)=\{(\sigma ,\beta ):\sigma\in\cJ,\beta\in\cM,
|\sigma|  +|\beta |\leq i\},
$$
$$
 B(i,j)=\{(\alpha ,\beta ):\alpha \in\cN ,\beta\in \cM  ,
|\alpha |=i,|\beta |\leq j\}.
$$

\begin{lemma}
                                                    \label{lemma 9.29.3}
 Let $\kappa,\mu\geq0$ be   
integers and 
$\alpha\in\cN,\beta\in \cM $,  
$\sigma\in\cJ$.   
Assume that 
\begin{equation}
                                                     \label{10.2.1}
2(|\sigma|+|\beta|+\kappa )+\mu +2 \leq l.
\end{equation}
Then
$$
S_{\sigma}(B_{\alpha}w_{\beta})=\sum_{i=0}^{\kappa  }
\delta^{i}\sumstar_{A(|\sigma|+|\beta|+i)} S_{\sigma_{1}}
w_{\beta_{1}}
$$
\begin{equation}
                                                \label{10.1.1}
+\sum_{ i=1}^{\kappa }\delta^{i}
\sumstar_{B(|\alpha|+i,|\sigma|+|\beta|+i-1)} 
B_{\alpha_{1}}w_{\beta_{1}}
+O_{\mu}( \delta^{\kappa+1}).
\end{equation}

\end{lemma}

\begin{proof} 
  For $\sigma=e$, when $S_{\sigma}=\cR$, equation (\ref{10.1.1}) 
turns out to be just a different form
of (\ref{9.29.2}), which is applicable since
$2(|\beta|+\kappa)+1\leq l$. Indeed,
owing to Lemma \ref{lemma 19.20.10.03} (iii) 
$$
\sum_{|\gamma|=i}c_{\alpha0\gamma}
 \cR 
w_{\beta\gamma}=
\sumstar_{A(|\sigma|+|\beta|+ i)} 
S_{\sigma_{1}}w_{\beta_{1}},
$$

$$
\sum_{|\gamma|=i-1}
B_{\alpha0\gamma} 
w_{ \beta\gamma}=
\sumstar_{B(|\alpha|+i,|\sigma|+|\beta|+ i-1)} 
B_{\alpha_{1}}w_{\beta_{1}}.
$$
Furthermore, for $|\gamma|=\kappa$ (see Remark \ref{remark 10.12.1})
$$
 \cQ_{\alpha0\gamma}
w_{ \beta\gamma}=O_{\mu}(1),\quad \text{since}
\quad  2(|\beta|+\kappa )+\mu+2\leq l.
$$

 For $|\sigma|\geq1$
we  proceed by induction on
the length $\ell(S_{\sigma})$ of $S_{\sigma}=
 \cR L_{\beta_{1}} 
\cdot...\cdot 
 \cR L_{\beta_{j}} $, 
which we define to be $j$. If $\ell(S_{\sigma})=1$, then
$S_{ \sigma }=
 \cR L_{\nu} $ for a $\nu\in \cM $
with $\nu=\sigma$ and it suffices to notice that
\begin{equation}
                                                     \label{10.29.1}
S_{ \sigma }
(B_{\alpha}w_{\beta})
 =\cR L_{\nu}(B_{\alpha}w_{\beta})
=-\cR(B_{\alpha}w_{\nu\beta})=
-S_{e}(B_{\alpha}w_{\beta'}),
\end{equation}
where $\beta'=\nu\beta\in\cM$ and $2(|\beta'|+\kappa)+\mu+2
=2(|\sigma|+|\beta|+\kappa)+\mu+2\leq l$.

Assume that (\ref{10.1.1}) holds 
whenever $\ell(S_{\sigma})=s$, and take an $S_{\sigma}$
such that $\ell(S_{\sigma})=s+1$. Then
$S_{\sigma}= \cR 
L_{\nu}  S_{\sigma'}$, where $\nu,\sigma'\in\cM$,
$|\nu|+ |\sigma'| =|\sigma|$ 
and $\ell(S_{\sigma'})=s$. Furthermore,
$$
2( |\sigma'| +|\beta|+\kappa)+\mu'+2\leq l,
$$
where $\mu'=\mu+2|\nu|$. By the induction hypothesis
$$
S_{\sigma'}(B_{\alpha}w_{\beta})=\sum_{i=0}^{\kappa}
\delta^{i}\sumstar_{A( |\sigma'|+|\beta|+i)}
 S_{\sigma_{1}}w_{\beta_{1}}
$$

$$
+\sum_{ i=1}^{\kappa }\delta^{i}
\sumstar_{B(|\alpha|+i, |\sigma'|+|\beta|+i-1)} 
B_{\alpha_{1}}w_{\beta_{1}}
+O_{\mu'}( \delta^{\kappa+1}).
$$
We apply 
$ \cR L_{\nu} $ 
to both parts of this
equality and take into account that
$L_{\nu}w_{\beta_1}=-w_{\nu\beta_1}$ and  
$|\nu|+ |\sigma'| =|\sigma|$. Then  similarly to
(\ref{10.29.1}) we get  that 
$$
S_{\sigma }(B_{\alpha}w_{\beta})
=\sum_{i=0}^{\kappa }
\delta^{i}\sumstar_{A( |\sigma|  +|\beta|+i)}
 S_{\sigma_{1}}w_{\beta_{1}}
$$
\begin{equation}
                                                 \label{10.2.2}
+\sum_{ i=1}^{\kappa }\delta^{i}
\sumstar_{B(|\alpha|+i,  |\sigma| +|\beta|+i-1)} 
 S_{e} (B_{\alpha_{1}} 
w_{\beta_{1}}) 
+O_{\mu}( \delta^{\kappa+1}).
\end{equation}
 
Now we transform the second term on the right.
Take $(\alpha_{1},\beta_{1})\in
B(|\alpha|+i,  |\sigma|  +|\beta|+i-1)$ and
notice that then 
$|\beta_{1}|\leq  |\sigma|  +|\beta|+i-1$. 
  Hence 
by assumption (\ref{10.2.1})  
$$
  2(|\beta_{1}|+\kappa-i)+\mu+2<l .
$$
Therefore, by the result for  $\sigma= e $
$$
S_{e} (B_{\alpha_{1} }
w_{ \beta_{1}})=
\sum_{j=0}^{\kappa-i}
\delta^{j}\sumstar_{A( |\beta_{1}|+j) }
 S_{\sigma_{2}}w_{\beta_{2}}
$$
$$
+\sum_{j=1}^{\kappa-i}\delta^{j}
\sumstar_{B(|\alpha_{1}|+j, |\beta_{1}|+j-1)} 
B_{\alpha_{2}}w_{\beta_{2}}
+O_{\mu}( \delta^{\kappa-i+1}).
$$
We substitute this result into (\ref{10.2.2})
and obtain (\ref{10.1.1}) after
collecting the coefficients of $\delta^{i+j}$ 
and noticing that,
if $(\alpha_{1},\beta_{1})\in
B(|\alpha|+i,  |\sigma|  +|\beta|+i-1)$ and 
$(\alpha_{2},\beta_{2})\in
B(|\alpha_{1}|+j, |\beta_{1}|+j-1 )$,
then
$$
(\alpha_{2},\beta_{2})\in B(|\alpha|+i+j,
|\sigma|+|\beta|+i+j-1).
$$
This justifies the induction and finishes 
the proof of the lemma.
\end{proof}

We remind the reader that
throughout this section 
the assumptions of Theorem \ref{theorem 11.26.06} are supposed to be
satisfied and 
in the following proposition  use the notation
$$
B^{*}(i,j)=\bigcup_{i_{1}=1}^{i}
B(i_{1},j).
$$

\begin{proposition}
                                             \label{lemma 9.30.1}
For any 
$j=0,1,...,k$ we have 
$$
w 
=v+\sum_{i=1}^{j}\delta^{i}
\sumstar_{ A(2i) }
 S_{\sigma} v_{ \beta }
+\sum_{i=j+1}^{ k }\delta^{i}
\sumstar_{A(i+j+1 )}S_{\sigma_{1}}
w_{\beta_{1}}
$$
\begin{equation}
                                              \label{10.2.4}
+\sum_{i=1}^{ k }
\delta^{i}\sumstar_{B^{*}(i,i+j)}
B_{\alpha_{1}}w_{\beta_{1}}
+O_{m}
(\delta^{ k +1}),
\end{equation} 
 where $v_{\beta}:=L_{\beta}v+f_{\beta}$.  
\end{proposition}

\begin{proof} 
By Proposition
\ref{01.24.06} (notice that, due to (\ref{10.4.5}),
$l\geq 2k+3$ and $2(k+1)+m+2
\leq l$)
we have
$$
w=v+\sum_{i=1}^{k}\delta^i
\sum_{|\beta|=i}B_{\beta}w_{\beta}
+\sum_{i=1}^{k}\delta^i\sum_{|\beta|=i+1}c_{\beta }
\cR w_{\beta}+O_{m}(\delta^{k+1}).  
$$
which means that
(\ref{10.2.4}) holds for $j=0$, since 
by Lemma  \ref{lemma 19.20.10.03} (iii) 
$$
\sum_{|\beta|=i}B_{\beta}w_{\beta}=
\sum_{|\beta|=i}B_{\beta}\sum_{\gamma\in\cM,|\gamma|=i}
c(\beta,\gamma)w_{\gamma}=\sumstar_{B^{*}(i,i)}
B_{\alpha_{1}}w_{\beta_{1}},
$$
$$
\sum_{|\beta|=i+1}c_{\beta }\cR w_{\beta}
=\sum_{|\beta|=i+1}c_{\beta }
\sum_{\gamma\in\cM,|\gamma|=i+1}c(\beta,\gamma)
\cR w_{\gamma}=
\sumstar_{A(i+1)}S_{\sigma_{1}}
w_{\beta_{1}}.
$$
  
Next, assume that 
$ k \geq1$ 
and (\ref{10.2.4}) holds for a
$j\in\{0,...,k-1\}$. 
Transform the first term
with
$i=j+1$ in the second sum on the right 
in (\ref{10.2.4}) by using Lemma
\ref{lemma 9.29.3}. 
To prepare the transformation 
take $(\sigma_{1},\beta_{1})
\in A(2i)=A(i+j+1)$
so that $|\sigma_{1}| 
+|\beta_{1}|\leq2i $ and 
apply the operator 
$S_{\sigma_{1}} L_{\beta_1} $ to both parts of
equation  
(\ref{02.24.06}) with  $k-i$ in place of 
$\kappa$.
Then
we obtain 
$$
S_{\sigma_{1}}w_{\beta_{1}}=
S_{\sigma_{1}}v_{\beta_{1}}
+\sum_{i_{1}=1}^{k-i}\delta^{i_{1}}
\sum_{|\alpha_{1}|=i_{1}}
S_{\sigma_{1}}(B_{\alpha_{1}}
L_{\beta_{1}}w_{\alpha_1})
$$
$$
+\sum_{i_{1}=1}^{k-i}\delta^{i_{1}}
\sum_{|\alpha_{1}|=i_{1}+1}c_{\alpha_{1}}
S_{\sigma_{1}}L_{\beta_{1}}\cR w_{\alpha_{1}}+
\delta^{k-i+1}r^{(k-i)},
$$
where
$$ 
r^{(k-i)} 
=\sum_{|\alpha|
=k-i+1}S_{\sigma_1}L_{\beta_1}\cQ_{\alpha}w_{\alpha}.
$$
Owing to 
$$
l-2(k-i+1
+|\beta_{1}|+ |\sigma_{1}| ) 
\geq l-
2(k+i+1 ) 
$$
$$
\geq l-2(2k
+1 )  \geq m+2, 
$$
we have 
 $
r^{(k-i)}= O_{m}(1)$.
 By the way, this is the only place where
we need $l$ to be not smaller than $4+m+4k$.
 Hence by Lemma  \ref{lemma 19.20.10.03} (iii)  
$$
 S_{\sigma_{1}}w_{\beta_{1}}= 
S_{\sigma_{1}}v_{\beta_{1}}
+\sum_{i_{1}=1}^{k-i}\delta^{i_{1}}
\sumstar_{(\alpha_{2},\beta_{2})\in B(i_{1},|\beta_{1}|+i_{1})}
S_{\sigma_{1}}(B_{\alpha_{2}}w_{\beta_{2}})
$$
\begin{equation}
                                                           \label{10.2.8}
+\sum_{i_{1}=1}^{k -i}\delta^{i_{1}}
\sumstar_{A( |\sigma_{1}| +|\beta_{1}|+i_{1}+1)} 
S_{\sigma_{2}}  w_{\beta_{2}}
+O_{m}
(\delta^{k-i+1})=:J_{1}+...+J_{4}.
\end{equation}
Now Lemma  \ref{lemma 9.29.3} with
$k-i-i_{1}$ in place of $\kappa$ and $m$ in place of $\mu$
allows  us to transform
 terms entering $J_{2}$. For $|\beta_{2}|
\leq|\beta_{1}|+i_{1}$ we have (remember that $
(\sigma_{1},\beta_{1})\in A(2i )$)
$$
2( |\sigma_{1}| +|\beta_{2}|+k-i-i_{1})+m
 +2 
\leq2(|\sigma_{1}| +|\beta_{1}|+k-i)+m +2 
$$
$$
\leq2(i+k )+m +2 
\leq4k+m +2 < l.
$$
Therefore  
$$
S_{\sigma_{1}}(B_{\alpha_{2}}w_{\beta_{2}})
=\sum_{i_{2}=0}^{k
-i-i_{1}}
\delta^{i_{2}}\sumstar_{A( |\sigma_{1}| +|\beta_{2}|+i_{2})}
 S_{\sigma_{3}}w_{\beta_{3}}
$$
$$
+\sum_{ i_{2}=1}^{k-i-i_{1}}\delta^{i_{2}}
\sumstar_{B(|\alpha_{2}|+i_{2}, |\sigma_{1}|+|\beta_{2}|+i_{2}-1)} 
B_{\alpha_{3}}w_{\beta_{3}}
+O_{m}( \delta^{k-i-i_{1}+1}).
$$
We plug this result into $J_{2}$ and in order to
collect the coefficients
of $\delta^{i_{1}+i_{2}}$ notice that,
for $(\sigma_{3},\beta_{3})\in A( |\sigma_{1}| +|\beta_{2}|
+i_{2})$ and $(\alpha_{2},\beta_{2})\in B(i_{1},|\beta_{1}|
+i_{1})$ it holds that
$$
 |\sigma_{3}| +|\beta_{3}|\leq |\sigma_{1}| +|\beta_{2}|+i_{2}
\leq |\sigma_{1}| +|\beta_{1}|+i_{1}+i_{2}.
$$
Furthermore, if $(\alpha_{3},\beta_{3})\in 
B(|\alpha_{2}|+i_{2}, |\sigma_{1}|+|\beta_{2}|+i_{2}-1)$, then
$$
|\alpha_{3}|=|\alpha_{2}|+i_{2}=i_{1}+i_{2},\quad
|\beta_{3}|\leq |\sigma_{1}| +|\beta_{2}|+i_{2}-1
< |\sigma_{1}|+|\beta_{1}|+i_{1}+i_{2}.
$$
It follows that $J_{2}$ is written as
$$
 \sum_{i_{1}=1}^{k-i }
\delta^{i_{1}}\big(
\sumstar_{A( |\sigma_{1}| +|\beta_{1}|+i_{1})}S_{\sigma_{2}}
w_{\beta_{2}}+\sumstar_{B(i_{1},|\sigma_{1}|+|\beta_{1}|+i_{1})}
B_{\alpha_{2}}w_{\beta_{2}}\big)+
O_{m}( \delta^{k-i+1}),
$$
which just amounts to saying that visually
  in the definition
of $J_{2}$ one can erase   $S_{\sigma_{1}}$, carry
all differentiations in it onto $w_{\beta_{2}}$,
and still preserve~(\ref{10.2.8}).

  Then we see that
$$
 \delta^{j+1}
\sumstar_{A(2j+2 )}S_{\sigma_{1}}w_{\beta_{1}}
=O_{m}( \delta^{k +1})
+\delta^{j+1}\sumstar_{A(2j+2 )}S_{\sigma_{1}}
v_{\beta_{1}}
$$
$$
+\sum_{i_{1}=1}^{k-j-1}
\delta^{i_{1}+j+1}
\big(
\sumstar_{A(|\sigma_{1}| 
+|\beta_{1}|+i_{1}+1)}S_{\sigma_{2}}
w_{\beta_{2}}
+\sumstar_{B(i_{1}, |\sigma_{1}| 
+|\beta_{1}|+i_{1})}
B_{\alpha_{2}}w_{\beta_{2}}\big).
$$
Next we notice again that, for $(\sigma_{1},\beta_{1})
\in A(2j+2 )$ and $ |\sigma_{2}| +|\beta_{2}|
\leq |\sigma_{1}|+|\beta_1|
 +i_{1}+1$, we
have
$ |\sigma_{2}|+|\beta_{2}|\leq  j+2 +i_{1}+j+1$,
whereas if $ |\beta_{2}|
\leq |\sigma_{1}| 
+|\beta_{1}|+i_{1}$, then
$ |\beta_{2}|\leq  j+1 +i_{1}+j+1$.
Therefore, after changing $i_{1}+j+1\to i$ 
($\geq j+2$) we get
$$
 \delta^{j+1}
\sumstar_{A(2j+2 )}
S_{\sigma_{1}}w_{\beta_{1}}
=O_{m}(\delta^{k+1})
+\delta^{j+1}\sumstar_{A(2j
+2 )}S_{\sigma_{1}}
v_{\beta_{1}}
$$
$$
+\sum_{i =j+2}^{k}
\delta^{i }
\big(
\sumstar_{A(i+j+2  )}S_{\sigma_{2}}
w_{\beta_{2}}+\sumstar_{B^{*}
(i ,i+j+1)}
B_{\alpha_{2}}w_{\beta_{2}}\big).
$$
This shows that the term with $i=j+1$ in 
the second sum on the right in (\ref{10.2.4})
can be eliminated on the account of changing
other terms with simultaneous shift $j\to j+1$.
Thus the induction on $j$ proves the 
proposition indeed.
\end{proof}

Now we finish the proof of 
Theorem \ref{theorem 11.26.06}. 
By taking $j=k$ 
in Proposition \ref{lemma 9.30.1}, we find
\begin{equation}
                                                         \label{11.19.1}
w= v+ \sum_{i=1}^{k}\delta^{i}
w^{(i)}
+\sum_{B^{*}(k,2k)} 
c(\alpha,\beta,\delta)B_{\alpha}w_{\beta}
+O_{m}(\delta^{k+1}),
\end{equation}
where 
$$
w^{(i)}:=\sumstar_{A(2i)}
S_{\sigma}v_{\beta}\in  C_{w}([0,T],W^{m}_{p}),\quad i=1,2,...,k,
$$
are independent
of $\delta$,  and $c(\alpha,\beta,\delta)$ 
are certain constants.
It is not hard to follow our computations 
in order to see that
$$
\sup_{t\in[0,d_{1}T]}\sup_{n,\delta=T/n }\delta^{-(k+1)}
\|O_{m}(\delta^{k+1})(t,\cdot)\|_{m,p}\leq N,
$$
where the constant $N$ depends only on $d,d_{1},T,K,k,m,p,l$. 
After that to finish the proof it only remains to recall that
$B_{\alpha}(jd_1\delta)=0$ 
for all integers $j\geq0$ and
$$
v(d_1t)=u(t), \quad w(d_1t)
=u_{n}(t)\quad\forall t\in
T_n =\{iT/n:i=0,1,2,...,n\}.
$$

\mysection{The case of time 
dependent coefficients}
                                               \label{section 18.01.09.03} 

We consider here the Cauchy problem 
(\ref{17.03.09.03})-(\ref{18.03.09.03})
for time dependent coefficients. 
We split, as before,  
the coefficients and the free terms into $d_1$ terms, 
$$
a^{ij}=\sum_{r=1}^{d_1}a^{ij}_r, \quad 
a^{i}=\sum_{r=1}^{d_1}a^{i}_r,\quad 
a=\sum_{r=1}^{d_1}a_r,\quad
f=\sum_{r=1}^{d_1}f_r, 
$$
define $\delta=T/n$, $t_{i}=t_i^n=\delta i$, $T_{n}=
\{t_{i}:i=0,1,...,n\}$,  
and   keep Assumptions 
\ref{assumption 10.09.01},  
\ref{assumption 10.11.09.03}. 
As before, we also denote
$
L_r:=a_r^{ij}D_{ij}+a^i_rD_i+a_{r}
$. 
 
One of splitting-up approximations 
$u_n(t)$ for $t\in T_n$ is defined 
as follows.
Let $\bS_{st}^{(r)}$ 
be the   operator 
mapping each function $\varphi$ of an appropriate class into the solution
of 
the problem 
$$
D _{t} 
v(t,x)=d_1L_rv(t,x)+d_1f_r(t,x), \quad t>s,\quad v(s,x)=\varphi(x).
$$ 
 Then the approximations are introduced according to 
$$
u_n(0):=u_0,
$$ 
\begin{equation}                                   \label{13.17.09.03}
u_n(t_{i+1})
:=\bS_{t_{i\bar d},t_{i+1}}^{(d_1)}
...\bS_{t_{i1},t_{i2}}^{(2)}
\bS_{t_i,t_{i1}}^{(1)}u_n(t_i), 
\qquad  i=0,1,2,...,n, 
\end{equation}
where $t_{ij}:=t_i+j\delta/d_1$, 
for $j =1,2,...,d_1-1$, $\bar d:=d_1-1$.

There are many other ways to extend 
the splitting-up approximations (\ref{switching}) 
to PDEs with time dependent data. 
Along with 
(\ref{13.17.09.03}) we also consider another 
approximation, which has the advantage 
that in each step 
we need to solve a time independent PDE, 
which is usually 
more convenient in practice than 
solving time dependent PDEs. 
This time we define the approximation $u_n$  
by 
$$
u_n(0):=u_0,
$$ 
\begin{equation}                                   \label{14.17.09.03}
u_n(t_{i+1}^n)
:=\bS_{\delta}^{(d_1)}(t_{i+1}^n)
...\bS_{\delta}^{(2)}(t_{i+1}^n)
\bS_{\delta}^{(1)}(t_{i+1}^n)u_n(t_{i}^n), 
\qquad i =0,1,2,...,n, 
\end{equation}
where  
$\bS_{\delta}^{(r)}(s)\varphi$ denotes 
the solution of the problem 
\begin{equation}                            \label{15.17.09.03}
D_{t}
v(t)=L_r(s)v(t)+f(s),
\quad t\geq 0,  \qquad                      
v(0)=\varphi, 
\end{equation}
with 
$$
L_r(s):=a^{ij}_r(s,x)D_{ij}
+a^{i}_r(s,x)D_{i}+a_r(s,x), 
$$
for $r=1,2,...,d_1$. 
Notice that the coefficients of the operator 
$L_r(s)$ and $f(s)$ are ``frozen" at time 
$s$, thus (\ref{15.17.09.03}) is 
a Cauchy problem with 
time independent data. 

We extend Theorem 
\ref{theorem 11.26.06} as follows. 

\begin{theorem}
                                                  \label{theorem 15.15.09.03}
Let $m\geq 0$ and $k\geq 0$ 
be any integers. Let 
Assumptions \ref{assumption 10.09.01}  and 
\ref{assumption 10.11.09.03}  
hold with 
$
l\geq 4+m+4k. 
$
Let the splitting-up approximation $u_n$ 
be defined by  (\ref{13.17.09.03}) or 
by (\ref{14.17.09.03}). 
Then there exist functions $u^{(j)}\in
C_{w}([0,T],W^m_p)$, $j=1,2,...,k$,
$R^{(k)}_n\in C_{w}([0,T],W^m_p)$, such that 
$$
u_n(t,x)=u(t,x)
+\delta u^{(1)}(t,x)
$$
\begin{equation}                                  \label{25.23.10.03}
+\delta^2u^{(2)}(t,x)+...
+\delta^k u^{(k)}(t,x)+R^{(k)}_n(t,x)                       
\end{equation}
for all $t\in T_n$, 
$x\in\bR^d$, and $n\geq 1$.
The functions
$u^{(j)}$, $j=1,2,...,k$, 
are independent of $n$, and
$$                           
\sup_{t\in[0,T]}
\|R^{(k)}_n(t)\|_{m,p}\leq N\delta^{k+1}
$$
for all $n$, where $N$ depends only on 
$k,d,d_1,K,m,p,T$.
\end{theorem}

Clearly Theorem \ref{theorem 15.15.09.03} 
implies that we can again accelerate the 
convergence of the splitting-up approximations 
by considering 
$$
v_n(t,\cdot):=\sum_{j=0}^kb_ju_{2^jn}(t,\cdot),
\quad t\in T_n, 
$$
where $u_{2^jn}$ for all $j=0,1,..,k$ are defined  
by either (\ref{13.17.09.03}) or  by (\ref{14.17.09.03}). 
\begin{theorem}                                       
                                               \label{theorem 12.17.09.03} 
Let $m\geq 0$ and
$k\geq 0$ be any integers.  
Let Assumptions 
\ref{assumption 10.09.01} 
and 
\ref{assumption 10.11.09.03}   
hold with $\, l\geq 4+m+4k\,$.  
Then for all $n\geq1$
$$                                  
\max_{t\in T_n}\|v_n(t)-u(t)\|_{m,p}
\leq N\delta^{k+1},   
$$
where $N$ is a constant, 
depending only on 
$k,d,d_1,K,m,p,T$. 
\end{theorem}

Hence by Sobolev's Theorem on embedding 
$W^m_p$ into $C^{s}$ we 
immediately get the following result. 
\begin{theorem}                                     \label{Sobolev}
Let $m\geq 0$ and 
$k\geq 0$ be any integers.  
Let Assumptions 
\ref{assumption 10.09.01} 
and 
\ref{assumption 10.11.09.03}   
hold with $\, l\geq 4+m+4k\,$. 
Let $s\geq0$ be an integer such that 
$m\geq s+d/p$. 
 Then for all $n\geq1$
$$                                  
\max_{t\in T_n}\sup_{x\in\bR^d}\sum_{|\rho|\leq s}
|D^{\rho}v_n(t,x)-D^{\rho}u(t,x)|
\leq N\delta^{k+1},   
$$
where $N$ is a constant, 
depending only on 
$k,d,d_1,K,m,s,p,T$. 
\end{theorem}

We prove Theorem \ref{theorem 15.15.09.03} 
by adapting the proof of Theorem \ref{theorem 11.26.06}
to the time dependent case. 
If $u_n$ is defined by (\ref{13.17.09.03}),  
then we consider the problems 

\begin{equation}                        \label{equation 16.16.09.03}
dv(t)=
\big(L(A_{0}(t))v(t)
+f(A_{0}(t))\big)\,dA_{0}(t),\quad
v(0)=u_0,
\end{equation} 
\begin{equation}                           \label{equation 15.16.09.03} 
dw(t)=
\sum_{r=1}^{d_1}\big(L_r(A_{0}(t))w(t)
+f_r(A_{0}(t))\big)\,dA_r(t),\quad
w(0)=u_0, 
\end{equation} 
where $ A_0(t)$,  
$A_1(t),A_2(t),...,A_{d_1}(t)$  are 
defined by  (\ref{9.24.1}) and (\ref{10.17.09.03}),  
and $L(A_{0}(t))$, $L_r(A_{0}(t))$ mean that we substitute 
$A_{0}(t)$ in place of the time variable 
$t$ of the coefficients 
of $L$, $L_r$. 

If $u_n$ is defined by (\ref{14.17.09.03}) 
then we consider problems 
(\ref{equation 16.16.09.03}) and 
(\ref{equation 15.16.09.03}) with 
absolutely continuous functions 
$ A_0$, $A_1$,...,$A_{d_1}$,  
defined by the following requirements: 
$$
A_r(0)=0,\quad \dot A_{r} 
{\text{  \rm is periodic with period }}  (d_1+1)\delta,   
$$
\begin{equation}                                 \label{01.24.14.03}                                      
\quad \dot A_{r}(t)=1_{[r,r+1]}(t/\delta), 
\quad t\in[0,(d_1+1)\delta] \quad{\text{ {\rm for }}} 
\quad  
r=0,1,...,d_1. 
\end{equation}

By virtue of Theorem \ref{theorem28.11.01} 
equations (\ref{equation 16.16.09.03})  
and (\ref{equation 15.16.09.03}) admit 
 unique solutions $v$ and $w$, respectively.   
Clearly  $v, w\in C_w([0,d^{'}T],W^l_p)$, and  
$$ 
v(d^{\prime}t)=u(t), \quad w(d^{\prime}t)=u_n(t)
\qquad \text{for all }t\in T_n,  
$$ 
where $d^{\prime}=d_1$ if $ A_0,A_1,...,A_{d_1}$ 
are defined by (\ref{9.24.1}) and (\ref{10.17.09.03}), 
and $d^{\prime}=d_1+1$ if $ A_0,A_1,...,A_{d_1}$ are defined 
by (\ref{01.24.14.03}).  
Therefore, our aim is to get an equality 
like (\ref{25.23.10.03})  
with $v$ and $w$ in place of 
$u$ and $u_n$, respectively.   

We treat the cases of two approximations simultaneously
and warn the reader that,
in order not to repeat the same arguments twice,
 we are going to use the same notation
for some objects that have different meaning in each case. 
From now on 
$d^{\prime}$ denotes $d_1$ if we consider 
the splitting-up approximations $u_n$ 
defined by (\ref{13.17.09.03}), and it denotes 
$d_1+1$ in the case of $u_n$ defined by (\ref{14.17.09.03}).
We keep the notation $\cN$ for the set of all multi-numbers 
$\alpha=\alpha_1\alpha_2...\alpha_j$ for  
$\alpha_i\in\{0,1,2,...,d_1\}$ and integers $j\geq1$.
 We use also 
the numbers $c_{\alpha}$ and  the functions $B_{\alpha}$, 
defined by (\ref{21.16.09.03}), (\ref{18.19.09.03}), and
(\ref{22.16.09.03}), with $d^{\prime}$ in place of 
$d_1$ in (\ref{18.19.09.03}). 
Observe that as is easy to see Lemma \ref{lemma10.25.06}
still holds with $d'$ in place of $d_{1}$ in its formulation.

Let $\cR f$ and $\bar \cR f$ 
denote the solutions of the problems 
$$
du(t)=(Lu(t)+f(t))\,dt,\quad u(0)=0, 
$$
and  
$$
dv(t)=(\bar Lv(t)+f(t))\,d A_{0}(t),\quad v(0)=0, 
$$
respectively, where $\bar L:=L(A_{0}(t))$, 
the operator obtained from $L$ by the substitution 
of $A_{0}(t)$ in place of $t$ in the coefficients 
of $L$. Notice that $\bar \cR$ depends on $n$ when 
$A_{0}$ is defined by (\ref{01.24.14.03}). Notice also 
that 
\begin{equation}                             \label{18.23.10.03}
\overline{\cR  f} (t,\cdot):=(\cR f)(A_{0}(t),\cdot)
=\bar \cR \bar f (t,\cdot), 
\end{equation}
where $\bar f (t,\cdot)=f(A_{0}(t),\cdot)$. 
Let $\bar Q_{\alpha}f$ denote the solution of the problem 
$$
dv(t)= \bar Lv(t)\,dA_{0}(t)
+f(t) \,dB_{\alpha}(t),\quad v(0)=0.  
$$

We modify the definition of $L_{\alpha}$, $f_{\alpha}$, used 
for time independent operators and free term, as follows: 
$L_0:=0$, $f_0:=0$,  
$$
L_{\gamma}:=L_r,\quad f_{\gamma}:=f_r
$$
for $\gamma=r\in\{1,2,...,d_1\}$, and 
$$
L_{\gamma0}:=LL_{\gamma}-\dot L_{\gamma},
\quad L_{\gamma r}:=-L_{\gamma}L_r
$$
\begin{equation}
                                                 \label{10.31.1}
f_{\gamma0}:=Lf_{\gamma}-\dot f_{\gamma},
\quad f_{\gamma r}:=-L_{\gamma}f_{r}
\end{equation}
for $r=1,2,..., d_1$, where 
$
\dot f_{\gamma}:=D_{t}f_{\gamma}, 
$ 
and $\dot L_{\gamma}$ denotes the differential operator which 
we obtain from $L_{\gamma}$ by taking the derivative in
$t$ of its coefficients. 
As is easy to see, $L_{\gamma}$ and $f_{\gamma}$
are well defined if $2(|\gamma|-1)\leq l$.
We use the notation 
$\bar L_{\gamma}$ and $\bar f_{\gamma}$
for the operator which we obtain from $L_{\gamma}$ 
by substituting 
$A_{0}(t)$ in place of $t$ in its coefficients, 
and for the function obtained from $f_{\gamma}$ 
by the same substitution, respectively. 
Then we have the following counterpart of 
Lemma  \ref{lemma 11.09.07.02}.
\begin{lemma}                                   \label{lemma 12.22.10.03}

Take some functions
$$
h \in L_{p}([0,d^{\prime}T], W^{1}_{p}),\quad
 h_{r}\in L_{p}([0,d^{\prime}T], W^{1}_{p}),
\quad r=0,1,...,d_{1}.
$$
 Let
$u$ be a solution of the ``equation''
$$                                                   
du=h_{r}\,dA_{r}  =\sum_{r=0}^{d_{1}}h_{r}\,dA_{r}
$$
with $u(0)\in W^{1}_{p}$. Assume that 
$Lu\in L_{p}([0,d^{\prime}T], W^{1}_{p})$.
Then  for any $\alpha\in\cN$ 
$$                                             
\bar\cR(B_{\alpha} h)=c_{\alpha0}\bar\cR h 
+\delta \bar\cQ_{\alpha0} h, 
$$
$$                                                    
\bar\cQ_{\alpha}u
=\bar\cR(c_{\alpha0} \bar Lu -c _{\alpha r}h_{r})
+\delta\bar{\cQ}_{\alpha0} \bar Lu 
-\delta\bar{\cQ}_{\alpha r} h_{r} +B_{\alpha}u.
$$
\end{lemma}

The proof of this lemma is an obvious modification 
of that of Lemma ~\ref{lemma 11.09.07.02}. 

Next, let us use the notation 
$$
w_{\beta}=\bar L_{\beta}w+\bar f_{\beta}. 
$$
Since $w\in C_w([0,d^{'}T],W^l_p)$,
the functions $w_{\beta}$ are well defined for
$2|\beta|\leq l$. Under the same condition
the coefficients of $L_{\beta}$ and $f_{\beta}$
have the first derivative in time and these 
derivatives are under control. Furthermore,
$
dw_{\beta}=h_{r}\,dA_{r}$, where, as long as $2|\beta|+3\leq l$,
the functions
$$
h_{0}=(\bar{L}\bar{L}_{\beta}-\bar{L}_{\beta0})w
+\bar{L}\bar{f}_{\beta}-\bar{f}_{\beta0},
$$
$$
h_{r}=\bar{L}(\bar{L}_{r}w+\bar{f}_{r}),\quad r=1,...,d_{1},
$$
are bounded $W^{1}_{p}$-valued functions on $[0,d'T]$.

Then in the same way as Lemma \ref{lemma 9.28.1}, 
Proposition \ref{01.24.06} and Lemma \ref{lemma 9.29.1} 
are obtained by the aid of Lemma \ref{lemma 11.09.07.02}, 
using Lemma \ref{lemma 12.22.10.03} we get their 
counterparts, formulated as follows.  
\begin{lemma}
                                          \label{lemma 17.22.10.03}
Let $\alpha,\beta\in\cN$.
If $2|\beta|+3\leq l$ then
$$                                                                                    
\bar\cQ_{\alpha}w_{\beta}
=c_{\alpha r}\bar\cR w_{\beta r}
+\delta\bar\cQ_{\alpha r}w_{\beta r}+B_{\alpha}w_{\beta}.
$$
\end{lemma}

\begin{proposition}                           \label{proposition 18.22.10.03}
Let $\kappa\geq 0$ be an integer 
and $l\geq2\kappa+3$.  
 Then                                         
\begin{equation}                               \label{23.22.10.03} 
w=v+\sum_{i=1}^{\kappa}\delta^i
\sum_{|\alpha|=i}B_{\alpha}w_{\alpha}
+\sum_{i=1}^{\kappa}\delta^i\sum_{|\alpha|=i+1}c_{\alpha }
\bar\cR w_{\alpha}
+\delta^{\kappa+1}\sum_{|\alpha|=\kappa+1}
\bar\cQ_{\alpha}w_{\alpha}.
\end{equation}
 
\end{proposition}

\begin{lemma}
                                                 \label{lemma 22.22.10.03}
If $\kappa\geq0$ is an integer and $\alpha,\beta\in\cN$
and $2(|\beta|+\kappa)+1\leq l$,
then
$$
\bar\cR(B_{\alpha}w_{\beta})=\sum_{i=0}^{\kappa }
\delta^{i}\sum_{|\gamma|=i}c_{\alpha0\gamma}\bar\cR w_{\beta\gamma}
$$
$$
+\sum_{i=1}^{\kappa }\delta^{i}\sum_{|\gamma|=i-1}
B_{\alpha0\gamma}w_{\beta\gamma}
+\delta^{\kappa+1}\sum_{|\gamma|=\kappa }\bar\cQ_{\alpha0\gamma}
w_{\beta\gamma},
$$
where for any multi-numbers $\mu,\nu$
$$
\sum_{|\gamma|=0}c_{\nu\gamma}\bar\cR w_{\mu\gamma}:=
c_{\nu}\bar\cR w_{\mu},\quad
\sum_{|\gamma|=0}B_{\nu\gamma}w_{\mu\gamma}:=
B_{\nu}w_{\mu},
$$
$$
\sum_{|\gamma|=0}\bar{\cQ}_{\nu\gamma}
w_{\mu\gamma}:=\bar{\cQ}_{\nu}
w_{\mu}.
$$
\end{lemma}
In order to iterate equation (\ref{23.22.10.03}) we 
introduce the following class of indices. 
We say that 
\begin{equation}                                \label{26.22.10.03}
\beta =\gamma^{\nu}
:=\gamma_1^{\nu_1}\gamma_2^{\nu_2}
...\gamma_j^{\nu_j} 
\end{equation}
is a {\em graded multi-number\/} of length 
$|\beta| :=j+ \nu_1 +{\nu_2}+...+\nu_j $, 
if $\gamma_i\in\{1,2,...,d_1\}$, 
$\nu_i\geq0$ is any integer for $i=1,2,...,j$, 
where $j\geq1$ is any integer. 
If $\nu_i=0$ for some $i$, then we 
also write $\gamma_i$ in place of 
$\gamma_i^{0}$ in 
(\ref{26.22.10.03}). Let $\cK$ denote the 
set of all graded multi-numbers. For each 
$\beta=\gamma^{\nu}
=\gamma_1^{\nu_1}\gamma_2^{\nu_2}
...\gamma_j^{\nu_j}\in\cK$ of 
length $|\beta|\leq 1+l/2$ we introduce 
the following operators and functions: 
$$
L_{\beta}=L_{{\gamma}^{ \nu }}
:=(-1)^{|\beta|-1}L_{\gamma_1}^{(\nu_1)}
L_{\gamma_2}^{(\nu_2)}\cdot...\cdot 
L_{\gamma_{j}}^{(\nu_j)},
$$ 
\begin{equation}                                  \label{24.22.10.03} 
f_{\beta}=f_{{\gamma}^{ \nu }}
:=(-1)^{|\beta|-1}L_{\gamma_1}^{(\nu_1)}\cdot...
\cdot L_{\gamma_{j-1}}^{(\nu_{j}-1)}
f_{\gamma_{j}}^{(\nu_j)}, 
\end{equation} 
where 
$f_r^{(s)}
:=D^{s}_{t}f_r$, 
and $L_r^{(s)}$ denotes the operator which we obtain from 
$L_r$ by applying the derivation 
$D^{s}_{t}$ to each of its coefficients. 
By definition $f_r^{(0)} =f_r$  
and $L_r^{(0)} =L_r$. It is easy to see
that for $\beta\in\cN$, when $\beta=\beta^{0}\in\cK$,
definitions (\ref{24.22.10.03}) are consistent with
(\ref{10.31.1}).

\begin{lemma}                                 \label{lemma 25.22.10.03}
The following 
statements hold.

(i) Let $\beta,\gamma\in\cK$  be such that 
$|\beta|+|\gamma|\leq1+l/2$. Then 
$$                             
L_{\beta}L_{\gamma}=-L_{\beta\gamma},
\quad
L_{\beta}f_{\gamma}=-f_{\beta\gamma}. 
$$

(ii) Let $\alpha\in\cN$  be 
such that $\rho:=|\alpha|\leq1+l/2$. 
Then there 
exist constants
$c(\gamma)=c(\alpha,\gamma)\in\{0,\pm1,\pm2,..\}$ 
defined for all $\gamma\in\cK$  with 
$|\gamma|=\rho$, such that 
\begin{equation}                                \label{equation 25.20.10.03}
L_{\alpha}
=\sum_{\gamma\in\cK,|\gamma|=\rho}c(\gamma)L_{\gamma}, 
\quad 
f_{\alpha}
=\sum_{\gamma\in\cK,|\gamma|=\rho}c(\gamma)f_{\gamma}. 
\end{equation}
\end{lemma}

\begin{proof}
Part (i) follows immediately from the definition 
(\ref{24.22.10.03}) 
of $L_{\beta}$, $f_{\beta}$. 
Part (ii) clearly holds for $\alpha=0$
and 
$\alpha=r\in\{1,...,d_1\}$. Assume that equations 
(\ref{equation 25.20.10.03}) hold for 
some $\alpha\in\cN$, $|\alpha|<1+l/2$. Then 
$$
L_{\alpha r}
=-L_{\alpha}L_r=-\sum_{\beta\in\cK,|\beta|=|\alpha|}
c(\beta)L_{\beta}L_r
=\sum_{\beta\in\cK,|\beta|=|\alpha|}
c(\beta)L_{\beta r},
$$
$$
f_{\alpha r}
=-L_{\alpha}f_r
=-\sum_{\beta\in\cK,|\beta|=|\alpha|}c(\beta)L_{\beta}f_r
=\sum_{\beta\in\cK,|\beta|=|\alpha|}c(\beta)f_{\beta r}
$$
for $r\in\{1,2...,d_1\}$, and 
$$
L_{\alpha0}=LL_{\alpha}-\dot L_{\alpha}=\sum_{r=1}^{d_1}
\sum_{\beta\in\cM,|\beta|=\rho}c(\beta)L_rL_{\beta}-
\sum_{\gamma^{\nu}\in\cM,|\gamma^{\nu}|=\rho}
c(\gamma^{\nu})\dot L_{\gamma^{\nu}}, 
$$
$$
f_{\alpha0}=Lf_{\alpha}-\dot f_{\alpha}=\sum_{r=1}^{d_1}
\sum_{\beta\in\cM,|\beta|=\rho}c(\beta)L_rf_{\beta}
-\sum_{\gamma^{\nu}\in\cM,|\gamma^{\nu}|=\rho}
c(\gamma^{\nu})\dot f_{\gamma^{\nu}}. 
$$
Hence by using assertion (i) and
noticing that
$$
\dot L_{\gamma^{\nu}}=\sum_{|\mu|=1}L_{\gamma^{\nu+\mu}}, \quad 
\dot f_{\gamma^{\nu}}=\sum_{|\mu|=1}f_{\gamma^{\nu+\mu}},
$$ 
we get equations (\ref{equation 25.20.10.03}) for $\alpha r$ 
with $r=0,1,...,d_1$. Thus the 
 induction on the length of $\alpha$ 
completes the proof. 
\end{proof}

For $\beta\in\cK$ we write $\bar L_{\beta}$ and 
$\bar f_{\beta}$, when the time change 
$A_{0}(t)$ is done in the coefficients of $L_{\beta}$ 
and in $f_{\beta}$. 
We set $w_{\gamma}:=\bar L_{\gamma}w+\bar f_{\gamma}$ 
for $\gamma\in\cK$, $|\gamma|\leq 1+l/2$.  Notice 
that Lemma \ref{lemma 25.22.10.03} has an obvious 
translation in terms of these functions. 
Namely, by Lemma \ref{lemma 25.22.10.03} (ii) 
for every $\alpha\in\cN$ such that $\rho:=|\alpha|\leq1+l/2$  
there exist constants
$c(\gamma)=c(\alpha,\gamma)\in\{0,\pm1,\pm2,..\}$ 
defined for all $\gamma\in\cK$  with 
$|\gamma|=\rho$, such that 
$$                              
w_{\alpha}
=\sum_{\gamma\in\cK,|\gamma|=\rho}c(\gamma)w_{\gamma}.  
$$

\noindent
For every integer 
$i\geq 1$ we introduce finite sequences 
$\sigma:=(\beta_1,\beta_2,...,\beta_i)$ of 
graded multi-numbers $\beta_i\in\cK$, and we set 
$|\sigma|:=|\beta_1|+|\beta_2|+...+|\beta_i|$. We also introduce the 
empty sequence $e$ of length $|e|:=0$. The set of all 
these sequences is denoted by $\cI$. 
For $\sigma=(\beta_1,\beta_2,...,\beta_i)$ with $|\sigma|\leq 1+l/2$  
we define 
$$
S_{ \sigma}:=\cR
L_{\beta_{1}}\cR L_{\beta_{2}}\cdot...\cdot
\cR L_{\beta_{i}}, 
\quad
\bar S_{ \sigma}:=\bar{\cR}
\bar L_{\beta_{1}}\bar\cR\bar L_{\beta_{2}}\cdot...\cdot
\bar\cR \bar L_{\beta_{i}}, 
$$
and for $|\sigma|=0$ we set 
$$
S_{e}:=\cR, \quad \bar S_{e}:=\bar\cR. 
$$
Notice that for any  
$g\in L_p([0,T],W^{2|\sigma|}_p)$  
\begin{equation}                                             \label{14.23.10.03}
\bar S_{\sigma}\bar g(t,\cdot)
=(S_{\sigma}g)(A_{0}(t),\cdot), 
\end{equation}
where $\bar g(t,\cdot):=g(A_{0}(t),\cdot)$.  
This follows from (\ref{18.23.10.03}) by induction 
on $|\sigma|$. In order to formulate the counterparts 
of Lemma \ref{lemma 9.29.3} and Proposition \ref{lemma 9.30.1} 
we use the following sets 
$$
A(i)=\{(\sigma ,\beta )
:\sigma\in \cI,\beta \in\cK, |\sigma|  +|\beta |\leq i\},
$$
$$
 B(i,j)=\{(\alpha ,\beta ):\alpha \in\cN,\beta\in\cK,
|\alpha |=i,|\beta |\leq j\},
$$
$$
B^{*}(i,j)=\bigcup_{i_{1}=1}^{i}
B(i_{1},j). 
$$ 
Remember that if  
$g_{\nu}$ is a collection of functions indexed by a parameter
$\nu$ taking values in a set $A$, then    
$
\sumstar_{\nu\in A}g_{\nu}
$
means any linear combination of $g_{\nu}$ with coefficient
independent of the argument of $g_{\nu}$ and of $\delta$.

\begin{lemma}                                  \label{lemma 21.23.10.03}
Let $\sigma\in \cI$, $\kappa,\mu\geq0$
be integers, and 
$\alpha\in\cN,\beta\in\cK$.  Assume that 
$$
2(|\sigma|+|\beta|+\kappa )+\mu +2 \leq l.  
$$
Then
$$
\bar S_{\sigma}(B_{\alpha}w_{\beta})=\sum_{i=0}^{\kappa  }
\delta^{i}\sumstar_{A(|\sigma|+|\beta|+i)} \bar S_{\sigma_{1}}
w_{\beta_{1}}
$$
$$
+\sum_{ i=1}^{\kappa }\delta^{i}
\sumstar_{B(|\alpha|+i,|\sigma|+|\beta|+i-1)} B_{\alpha_{1}}w_{\beta_{1}}
+O_{\mu}( \delta^{\kappa+1}).
$$
\end{lemma}
\begin{proof}
We can derive this lemma from Lemma 
\ref{lemma 22.22.10.03} in the same way as Lemma 
\ref{lemma 9.29.3} is proved. We need only use the 
sets $\cK$ and $\cI$ in place of $\cM$ and $\cJ$, 
and the operators $\bar \cR$, $\bar L_{\nu}$, $\bar S_{\sigma}$
for 
$\nu\in\cK$, $\sigma\in \cI$,  
 in place of 
$\cR$, $L_{\nu}$ $S_{\sigma}$, for 
$\nu\in\cM$, $\sigma\in \cJ$, respectively.   
\end{proof}

\begin{proposition}
                                     \label{proposition 19.23.10.03}
Let 
$k, m\geq0 $ be integers,  
and 
$$                                              
4+m+4k \leq l.
$$
Then for any 
$j=0,1,...,k$ we have 
$$
w=v+\sum_{i=1}^{j}\delta^{i}
\sumstar_{A(2i)}\bar S_{\sigma}v_{\beta}
+\sum_{i=j+1}^{k}\delta^{i}
\sumstar_{A(i+j+1)}\bar S_{\sigma_{1}}
w_{\beta_{1}}
$$
$$
+\sum_{i=1}^{k}
\delta^{i}\sumstar_{B^{*}(i,i+j)}
B_{\alpha_{1}}w_{\beta_{1}}
+O_{m}(\delta^{k+1}),
$$
where $v_{\beta}:=\bar L_{\beta}v+\bar f_{\beta}$. 
\end{proposition}
\begin{proof}
The proof of this proposition is 
a straightforward translation of the proof of 
the corresponding proposition, 
Proposition \ref{lemma 9.30.1} 
in the time independent case. 
To make this translation we use the 
sets $\cK$ and $\cI$ in place of $\cM$ and $\cJ$, 
and the operators $\bar \cR$, $\bar L_{\nu}$, $\bar S_{\sigma}$
for 
$\nu\in\cK$, $\sigma\in \cI$,  
in place of 
$\cR$, $L_{\nu}$ $S_{\sigma}$, for 
$\nu\in\cM$, $\sigma\in \cJ$, respectively.   
\end{proof}

Now we can finish the proof of Theorem  
\ref{theorem 15.15.09.03} as follows. Taking 
$j=k$ in Proposition \ref{proposition 19.23.10.03},  
we get 
\begin{equation}                            \label{equation 24.23.10.03}
w=v+\sum_{i=1}^{j}\delta^{i}
\sum_{A(2i)}c(\sigma,\beta)\bar S_{\sigma}v_{\beta}                                        
+\sum_{B^{*}(k,2k)}c(\alpha,\beta,\delta)
B_{\alpha_{1}}w_{\beta_{1}}
+r_{\delta}, 
\end{equation} 
where $c(\sigma,\beta)$, $c(\alpha,\beta,\delta)$  
are certain constants, and 
$r_{\delta}$ is a function in $C_{w}([0,T],W^m_p)$
for each $\delta$, and 
$$
\sup_{t\in[0,d^{\prime}T]}
\sup_{n,\delta=T/n }\delta^{-(k+1)}
\|r_{\delta}(t,\cdot)\|_{m,p}\leq N. 
$$
Observe that in contrast with (\ref{11.19.1})
the functions $v$ and $\bar{S}_{\sigma}v_{\beta}$
 in (\ref{equation 24.23.10.03}) may depend
on $\delta$.
To proceed further, 
define $R^{(k)}_n(t,x):=r_{\delta}(d^{\prime}t,x)$, and   
$$
u^{(i)}
:=\sum_{A(2i)}c(\sigma,\beta)S_{\sigma}u_{\beta}, 
\quad i=1,2,..,k, 
$$
where $u_{\beta}:=L_{\beta}u+f_{\beta}$. 
Then by virtue of equality (\ref{14.23.10.03}) 
and the fact that $v(t)=u(A_{0}(t))$ from 
equation (\ref{equation 24.23.10.03}) we get 
$$
w(t,\cdot)=u(A_{0}(t),\cdot)+\sum_{i=1}^{j}\delta^{i}
u^{(i)}(A_{0}(t),\cdot)
$$
$$                                       
+\sum_{B^{*}(k,2k)}c(\alpha,\beta,\delta)
B_{\alpha_{1}}(t)w_{\beta_{1}}(t,\cdot)
+R^{(k)}_n(t/d ^{ \prime },\cdot). 
$$
Substituting here $d^{\prime}t$ in place of $t$ 
we get the required representation 
(\ref{25.23.10.03}) by taking into account that
$$
w(d^{\prime} t)=u_n(t),\quad A_{0}(d^{\prime}t)=t, 
\quad 
B_{\alpha_{1}}(d^{\prime}t)=0\quad 
\forall t\in T_n. 
$$

\begin{remark} Let $1\leq j\leq d_1$. 
Then Theorems \ref{theorem 15.15.09.03},  
\ref{theorem 12.17.09.03}, and 
\ref{Sobolev} hold also when the operator  
$\bS_{\delta}^{(r)}(t_{i+1}^n)$ is 
replaced with $\bS_{\delta}^{(r)}(t_{i}^n)$ 
for every $r=1,2,...,j$ 
in the definition (\ref{14.17.09.03}) of 
the splitting-up approximation $u_n$.  
To see this we need only repeat the proof 
of the previous theorem with  
 $A_j$ in place of $A_0$ 
in equation  (\ref{equation 15.16.09.03}) and 
with $A_{0}$ and $A_{j}$ interchanged in 
(\ref{equation 16.16.09.03}).   
\end{remark}

\end{document}